\documentstyle [11pt,graphicx,amstex,epsfig,amscd,amssymb]{article}
\def\al{\alpha}
\def\b{\beta}
\def\om{\omega}
\def\de{\delta}

\def\la{\lambda}

\def\phi{\varphi}
\def\ro{\rho}

\def\dim{\mbox{dim }}

\def\ti{\times}

\def\hf{{\mathfrak h}}
\def\sof{{\mathfrak so}}

\def\be{\begin{equation}}
\def\ee{\end{equation}}
\def\bear{\begin{eqnarray}}
\def\eear{\end{eqnarray}}
\def\best{\begin{eqnarray*}}
\def\eest{\end{eqnarray*}}
\def\pf{{\bf Proof}: }

\newtheorem{theorem}{Theorem}[section]

\newtheorem{proposition}[theorem]{Proposition}
\newtheorem{lemma}[theorem]{Lemma}
\newtheorem{cor}[theorem]{Corollary}

\newtheorem{defn}[theorem]{Definition}

\newtheorem{remark}[theorem]{Remark}

\newtheorem{example}[theorem]{Example}

\renewcommand{\Re}{\operatorname{Re}}
\renewcommand{\dim}{\operatorname{dim}}
\renewcommand{\Im}{\operatorname{Im}}
\def\non{\noindent}
\def\pf{\non {\bf Proof: }}

\def\ra{\rightarrow}

\def\r#1{\right#1}
\def\l#1{\left#1}
\def\ma#1{\mathop {#1} \limits}
\def\Z{{ \Bbb Z}}
\def\R{{ \Bbb R}}
\def\C{{ \Bbb C}}

\def\Q{{ \Bbb Q}}

\def\H{{ \Bbb H}}
\def\E{{ \Bbb E}}
\def\T{{ \Bbb T}}
\def\O{{ \Bbb O}}
\def\I{{ \Bbb I}}

\def\jbar{\bar\jmath}
\def\so{SO(4)}

\def\Scal{\mathcal S}
\def\a{\alpha}
\def\b{\beta}
\title{\bf Second Order Families of Special Lagrangian Submanifolds in ${\Bbb
C}^4 $}

\author{ Marianty Ionel\\McMaster University\\Hamilton, ON L8S 4K1}

\date{}
\begin{document}

%\tableofcontents
%\newpage

\maketitle

\vskip.15in
\begin{abstract}

This paper extends to dimension 4 the results in the article
 ``Second Order Families of Special
Lagrangian 3-folds'' by Robert Bryant. We consider the problem of
classifying the special Lagrangian 4-folds in $ \C^4 $ whose
fundamental cubic at each point has a nontrivial stabilizer in
$SO(4)$. Points on special Lagrangian 4-folds where the
SO(4)-stabilizer is nontrivial are the analogs of the umbilical
points in the classical theory of surfaces.

 In proving existence for the families of special Lagrangian 4-folds, we used the method
 of exterior differential systems in Cartan-K{\"a}hler theory. This method is guaranteed to
  tell us whether there are any families of
 special Lagrangian submanifolds with a certain symmetry,
but does not give us an explicit description of the submanifolds.
To derive an explicit description, we looked at foliations by
submanifolds and at other geometric particularities. In this
manner, we settled many of the cases and described the families of
special Lagrangian submanifolds in an explicit way.
\end{abstract}

%\pagebreak
 \tableofcontents

%%%%%%%%%%%%%%%%%%%%%%%%%%%%%%%%%%%%%%%%%%%%%%%%%%%%%%
%%%%%%%%%%%%% Introduction
%%%%%%%%%%%%%%%%%%%%%%%%%%%%%%%%%%%%%%%%%%%%%%%%%%%%%%
\pagebreak
\bigskip
\section{Introduction}
\bigskip
The complex space ${\C^m}$ is endowed with a K{\"a}hler form $$\om
 =\frac{i}{2}(dz_1\wedge d\bar z_1+dz_2\wedge d\bar z_2 +\dots
 +dz_m\wedge d\bar z_m) $$ and a volume form $$\Omega =
 dz_1\wedge dz_2\wedge \dots \wedge dz_m,$$ where $(z_1,z_2,\dots, z_m)$ are
 the coordinates on $\C^m$. A {\it special Lagrangian submanifold} in $ {\C^m}$
 is an $m$-dimensional real submanifold on which the forms $\om $ and
 $\Im\Omega$ restrict to 0.

The study of special Lagrangian submanifolds started with Harvey
and Lawson in their paper \cite{hl} on calibrated geometries. They
constructed many interesting examples of SL\footnotemark
\;$m$-folds in $\C^m $ and proved local existence theorems. Since
then, many other examples have been constructed using a variety of
techniques.  To give some
 examples, Dominic Joyce used the method of ruled submanifolds,
 integrable systems and evolution of quadrics in \cite{dj1}, \cite{dj2}, \cite{dj3} to construct
 explicit examples of special Lagrangian $m$-folds in $\C^m$, Mark Haskins exhibited
 examples of
 special Lagrangian cones in $\C^3$ \cite{ha}, Richard Schoen and Jon Wolfson used the
 variational approach for some of their constructions in Calabi-Yau manifolds in \cite{sw}, etc.

\footnotetext{Special Lagrangian}

Special Lagrangian geometry received reinforced attention in 1996
when Strominger, Yau and Zaslow formulated what is today known as
the SYZ conjecture \cite{sy}. This conjecture reveals the role of
the special Lagrangian geometry in Mirror Symmetry, a mysterious
relationship between pairs of Calabi-Yau 3-folds, coming from
String Theory. In this larger context, a lot of research is going
on nowadays to find examples of special Lagrangian submanifolds.
This would help in understanding what kind of singularities a
special Lagrangian submanifold in a Calabi-Yau can have,
classifying them and maybe ultimately resolving the SYZ
conjecture.

While, from the String Theory point of view, the most interesting
case to study is the special Lagrangian 3-folds of a Calabi-Yau
3-fold, higher dimensional cases are also important for the
understanding of the general theory of SL submanifolds in
Calabi-Yau $m$-folds.

The idea in this research, initiated by Robert Bryant in his paper
\cite{br1}, is to classify families of SL submanifolds that are
characterized by invariant, geometric conditions. When the ambient
space is flat, the second fundamental form is the lowest order
invariant of a SL submanifold, so we would like to study the
second order families of SL $m$-folds in $ \C^m $, that is the
families of SL $m$-folds in $ \C^m $ whose second fundamental form
satisfies a set of pointwise conditions.

 The second
fundamental form of a special Lagrangian submanifold in $ \C^m $
has a natural interpretation as a harmonic cubic form on the
submanifold, called the {\it fundamental cubic}. The stabilizer at
a generic point of the fundamental cubic of a generic SL $m$-fold
is trivial. For comparison, in the case of a hypersurface in $
\R^{m+1} $, the stabilizer of the second fundamental form in $
SO(m)$ is always nontrivial and is larger than the minimum
possible stabilizer exactly at the umbilical points of the
hypersurface. For this reason, the points on SL $m$-folds where
the $SO(m)$-stabilizer is nontrivial are the analogs of the
umbilical points in the classical theory of surfaces.

In his article \cite{br1}, Robert Bryant considered the
`umbilical' case and completely classified the SL submanifolds of
$ \C^3 $ whose fundamental cubic has nontrivial $
SO(3)$-stabilizer at a generic point. He found that the only SL
3-folds whose fundamental cubic has a nontrivial stabilizer at a
generic point are the 3-planes, with stabilizer $SO(3)$, the
Harvey and Lawson examples, with stabilizer SO(2), the austere SL
3-folds, with stabilizer $ S_3 $, the asymptotically conical SL
3-folds, with stabilizer $ \Z_3 $ and the Lawlor-Harvey-Joyce
examples, with stabilizer $ \Z_2 $.

This present work extends these results to dimension $ m=4 $,
namely tries to classify the special Lagrangian 4-folds in $ \C^4
$ whose fundamental cubic at a generic point has nontrivial $
SO(4)$-stabilizer.

The possible stabilizer of a harmonic cubic can be a continuous,
meaning a positive dimensional, or a discrete subgroup of $SO(4)$.
In chapter 3.2, we consider the case when the stabilizer is
continuous. It turns out that there are four cases when
 there are nontrivial special Lagrangian
 4-folds with continuous symmetry: (a) when the fundamental cubic has symmetry $SO(3)$ we obtain the
 Harvey and Lawson examples which appeared also in dimension 3: $L_c=\{(s+it){\bf u}| {\bf u}\in S^3\subset
 \R^4,\; \Im(s+it)^4=c\}$, where $c$ is any real constant; (b) when the symmetry
  is $SO(2)\ltimes S_3$, we obtain special Lagrangian submanifolds as products of the form $
 L=\R^2\ti \Sigma$, where $\Sigma\subset \C^2$ is a complex curve; (c) when the symmetry
 is $SO(2)$, we obtain the $SO(3)$-invariant special Lagrangian 4-folds and (d) when the symmetry is
 an $O(2)$, we obtain a two parameter family of solutions which we have not been able to
 integrate completely yet.

In chapter 3.3, we consider the case when the stabilizer of the
fundamental cubic is a discrete subgroup of $\so$. In chapter
3.3.1., we classify the SL 4-folds with polyhedral symmetry. It
turns out that the polyhedral subgroups of $\so$ that stabilize a
harmonic cubic in 4 variables are the tetrahedral subgroup $\T$,
the irreducibly acting octahedral subgroup $\O^+$ and the
irreducibly acting icosahedral subgroup $\I^+$. We show that the
special Lagrangian 4-folds whose stabilizer of its fundamental
cubic is isomorphic to the tetrahedral subgroup are the
Harvey-Lawson examples invariant under a torus action, the ones
whose stabilizer at a generic point is isomorphic to $\O^+$ are
the cones on flat 3-dimensional tori in the 7-sphere and that
there are no nontrivial special Lagrangian 4-folds whose
stabilizer of its fundamental cubic at a generic point is
isomorphic to $\I^+$.

Using the classification of the discrete subgroups of $\so$ from
chapter 3.1., it remains to analyze the cases when the cubic
stabilizer is a cyclic or a dihedral subgroup of $\so$. We show
that the discrete stabilizer can only have elements of order less
or equal to 6. Further, we show that if the stabilizer is discrete
and contains an element of order 6, 5 or 4, then there are no
special Lagrangian 4-folds in $\C^4$ with a cyclic or dihedral
symmetry.

When the stabilizer contains an element of order 3, there are two
inequivalent orbits in the space of fixed harmonic cubics that
have to be considered. In the first case of discrete symmetry at
least a $\Z_3$, the special Lagrangian 4-folds whose cubic
stabilizer at a generic point is isomorphic to $D_3$, the dihedral
group in 3 elements, turn out to be asymptotically conical. The SL
4-folds with symmetry an order 18 normal subgroup of $D_3\ti D_3$
turn out to be products of holomorphic curves. When the symmetry
is exactly a $\Z_3$, we were able to show that there is an
infinite parameter family of solutions that depends on 4 functions
of 1 variable, foliated by minimal Legendrian surfaces and by
holomorphic curves, but could not finalize the analysis and
describe this family completely.

In the second case of discrete symmetry at least $\Z_3$, we found
a large family of SL 4-folds defined by holomorphic differential
equations, family that did not appear in dimension 3.

The general case of discrete symmetry at least a $\Z_2$ is the
most complicated case, since the general cubic has a large number
of parameters, and was not considered in this work.

{\it \bf Acknowledgements: } I would like to thank my advisor
Prof. Robert Bryant for introducing me to the subject, for his
help and support and for the innumerable hours of discussions that
led to this work.

%%%%%%%%%%%%%%%%%%%%%%%%%%%%%%%%%%%%%%%%%%
%%%%%%Special Lagrangian Geometry
%%%%%%%%%%%%%%%%%%%%%%%%%%%%%%%%%%%%%%%%%
\section{Special Lagrangian Geometry in Calabi-Yau manifolds}
\subsection{Special Lagrangian Submanifolds}
We begin with the definition of a Calabi-Yau manifold.

\begin{defn}
A Calabi-Yau m-fold $(M,J,g)$ is a compact, complex
$m$-dimensional manifold $(M,J)$ with trivial canonical bundle
$K_M$ and Ricci-flat K{\"a}hler metric $g$.
\end{defn}
Because the canonical bundle $K_M$ is trivial, there is a nonzero
holomorphic section $\Omega$ of $K_M$. Since the metric $g$ is
Ricci-flat, $\Omega$ is a parallel tensor with respect to the
Levi-Civita connection $\nabla^g$ \cite{dj}. By rescaling
$\Omega$, we can take it to be the holomorphic $(m,0)$-form that
satisfies:
\begin{equation}
{\frac{\om^m}{m!}}={(-1)^{\frac{m(m-1)}{2}}{\biggl(\frac{i}{2}\biggl)}^m\Omega\wedge\bar
\Omega}, \label{holo}
\end{equation}
where $\om$ is the K{\"a}hler form of $g$. The form $\Omega$ is
called the {\it holomorphic volume form} of the Calabi-Yau
manifold $M$.

The special Lagrangian submanifolds were introduced by Harvey and
Lawson in their paper \cite{hl} using calibrations. They are
defined in the general setting of a Calabi-Yau manifold and are a
special class of minimal submanifolds.

\begin{defn}
Let $(M,J,g,\Omega)$ be a Calabi-Yau m-fold and $L\subset M$ a
real m-dimensional submanifold of $M$. Then $L$ is called a {\it
special Lagrangian submanifold} of $M$ if $ \om\mid_L\equiv 0$ and
$\Im\Omega\mid_L\equiv 0$.
\end{defn}

More generally, $L$ is said to be a special Lagrangian submanifold
with phase $e^{i\theta}$ if $ \om\mid_L\equiv 0$ and
$\Im(e^{i\theta}\Omega)\mid_L\equiv 0$.

 As an example, we can see that $\R^m\subset\C^m$ is a
special Lagrangian subspace. $\C^m$ is endowed with the standard
Calabi-Yau structure defined by
\begin{align}
\label{str1}
g_0 & =dz_1\circ d\bar z_1+\dots + dz_m\circ d\bar z_m \\
\label{str2} \om_0 & =\frac{i}{2}(dz_1\wedge d\bar z_1+\dots
+dz_m\wedge d\bar
z_m)\\
\label{str3} \Omega_0 & = dz_1\wedge\dots\wedge dz_m
\end{align}
where $g_0$ is the K{\"a}hler metric on $\C^m$, $\omega_0$ the
K{\"a}hler form and $\Omega_0$ the holomorphic volume form on
$\C^m$. An $m$-dimensional submanifold $L$ of $M$ is called {\it
Lagrangian} if $\om\mid_L=0$. So, the special Lagrangian
submanifolds of $M$ are the Lagrangian submanifolds with the extra
condition $\Im\Omega\mid_L=0$, which is exactly the reason for
their name.

There are some important results on SL submanifolds which we will
briefly recall here.
\smallskip

\noindent {\bf a. Deformations.} R. McLean \cite{ml} studied the
moduli space of compact special Lagrangian deformations
 and showed that it has the following description.
\begin{theorem}
[McLean] Let $(M,J,g,\Omega)$ be a Calabi-Yau m-fold and $L\subset
M$ a m-dimensional compact SL submanifold. Then the moduli space
${\mathcal M}_L$ of special Lagrangian deformations of $L$ is a
smooth manifold of dimension $b^1(L)$, the first Betti number.
\end{theorem}
\noindent{\bf b. Local Existence.} Harvey and Lawson \cite{hl}
proved local existence only for SL-submanifolds in $\C^m$, but
their result extends to show that if $(M,J,g,\Omega)$ is a
Calabi-Yau $m$-fold and $N\subset M$ a real analytic submanifold
of dimension $m-1$ such that $i^*(\om)=0$, then $N$ lies in a
unique irreducible SL submanifold $L\subset M$. Here $i\colon N\to
M $ is the inclusion map.

This result shows that there are many special Lagrangian
submanifolds locally.
\smallskip

\noindent{\bf c. Minimizing Property.} A closed special Lagrangian
submanifold is volume-minimizing in
 its homology class and therefore it is a
minimal submanifold. We remark that a minimal submanifold, i.e. a
submanifold with
 constant mean curvature 0, is not necessarily volume-minimizing
 amongst homologous submanifolds. For
 example the equator of a  2-dimensional sphere is minimal, but does
 not minimize length amongst lines of latitude.
\bigskip
\subsection{Structure equations} {\bf a. The Coframe Bundle}

Let $(M,J,g,\Omega)$ be a Calabi-Yau $m$-fold and let
${\C^m}\cong{\R^{2m}}$ have complex coordinates
$(z_1,z_2,\dots,z_m)$ and complex structure $I$. The standard
Calabi-Yau structure on ${\C^m}$ is given by the relations
(\ref{str1}), (\ref{str2}) and (\ref{str3}). Let $\pi\colon P\to M
$ denote the bundle of $\C^m$-valued Calabi-Yau coframes, i.e. an
element of $P_x=\pi^{-1}(x)$ is a complex linear vector space
isomorphism $u\colon T_x M \to \C^m$ that satisfies
$\om_x=u^*(\om_0) $ and $\Omega_x=u^*(\Omega_0)$. Then $\pi\colon
P\to M $ is a principal right $SU(m)$-bundle over $M$ and the
right action is given by $R_a(u)=a^{-1}\circ u$ for $a
 \in SU(m)$. $P_x$ is the fiber at $x$ of the
Calabi-Yau coframe bundle $P$.

The canonical form $\xi$ is defined on the Calabi-Yau coframe
bundle $P$ by
$$\xi_u=u\circ(d\pi)_u\colon T_u P \to \C^m \quad \mbox{for } u \in P$$ where $
(d\pi)_u\colon T_u P \to T_{\pi(u)} M$ is the differential of
$\pi$ at $u$. The 1-form $\xi$ is $\C^m$-valued and we denote its
components by $\xi_i$, $i=1\dots m$. Then, on the bundle $P$ the
following equations hold:
\begin{equation}
\pi^*(\om)=\frac{i}{2}(\xi_1\wedge\bar \xi_1+\dots +\xi_m\wedge
\bar \xi_m) \quad \mbox{and  } \pi^*(\Omega)=\xi_1\wedge \dots
\wedge \xi_m
\end{equation}

By regarding the forms on $M$ embedded into the forms on $P$ via
the pullback, we can ignore $\pi^*$ in the above equations.

We define also the functions $e_i\colon P \to TM $ such that
$\xi_i(e_j)=\de_{ij}$. So, if $v\in T_u P$ then: $\
(d\pi)_u(v)=e_i(u)\xi_i(v)$. Cartan's first structure equation:
\begin{align}
\label{struc1} d\xi_i=-\psi_{i\jbar}\wedge\xi_j
\end{align}
defines $(\psi_{i\jbar})=\psi$, the $\frak{su}(m)$-valued 1-form
on $P$ called the {\it connection form}. In the flat case
$M=\C^m$, Cartan's second structure equation satisfied by the
connection form $\psi$ is:
\begin{align}
\label{struc2} d\psi=-\psi\wedge\psi
\end{align}
\vskip.5cm \noindent{\bf b. Special Lagrangian Submanifolds in
$\C^m$}

In this paper we are interested in special Lagrangian submanifolds
of $\C^m$, therefore we are considering only the flat case from
now on. When $M=\C^m$ with the standard Calabi-Yau structure
$(\C^m,J,g_0,\Omega_0)$, we denote the Calabi-Yau coframe bundle
by $\mbox{x}\colon P \to \C^m$ and regard the functions $e_i$ as
vector-valued functions on $P\cong\C^m \times SU(m)$. Then the
relations:
\begin{align}
\label{1}
d\mbox{x} =e_i\xi_i\\
\label{2} de_i =e_j\psi_{j\bar\imath}
\end{align}
give the 1-forms $\{\xi_i,\psi_{i\jbar}\}$ which form a basis for
the space of 1-forms on the frame bundle $P$.

To study the SL submanifolds of $\C^m$, we separate the two
structure equations (\ref{struc1}) and (\ref{struc2}) into real
and imaginary parts. We set $\xi_i=\om_i+i\eta_i$ and
$\psi_{i\jbar }=\al_{ij}+i\b_{ij}$. The first structure equation
(\ref{struc1}) becomes
\begin{equation}
d\om_i=-\al_{ij}\wedge\om_j+\b_{ij}\wedge\eta_j \quad
\mbox{and}\quad d\eta_i=-\b_{ij}\wedge\om_j-\al_{ij}\wedge\eta_j
\label{sstruc1}
\end{equation}
where we used Einstein's convention to sum over repeated indices.
Since $\psi$ is skew-hermitian with trace 0, it follows that
$\al=(\al_{ij})$ is skew-symmetric and $\b=(\b_{ij})$ is symmetric
with vanishing trace.

 When split into its real and imaginary parts, the
second structure equation (\ref{struc2}) becomes:
\begin{align}
\label{Gauss}
d\al_{ij}=-\al_{ik}\wedge\al_{kj}+\b_{ik}\wedge\b_{kj}\\
\label{Codazzi}
d\b_{ij}=-\b_{ik}\wedge\al_{kj}-\al_{ik}\wedge\b_{kj}
\end{align}

 Let $L\subset M$ be a special Lagrangian
submanifold. We are going to consider the bundle $P_L$ of
L-adapted coframes over $L$. This is defined as follows. Let $x\in
L$. A Calabi-Yau coframe at $x$, $u\colon T_x M \to \C^m$ is said
to be {\it L-adapted }if $u(T_x L)=\R^m\subset\C^m$ and $u\colon
T_x L \to \R^m$ preserves orientation. The space of $L$-adapted
coframes forms a principal right $SO(m)$-subbundle
$P_L\subset\pi^{-1}(L)\subset P$ over $L$. Now, because $u$ takes
a tangent plane to $L$ in $M$ into a real one, $\xi$ is
$\R^m$-valued on $P_L$ and so $\eta_i=0$ holds on $P_L$. By the
structure equation (\ref{sstruc1}), we get that:
\begin{equation*}
d\om_i=-\al_{ij}\wedge\om_j \quad \mbox{and} \quad
\b_{ij}\wedge\om_j=0 \quad \mbox{on } P_L
\end{equation*}
 Since ${\om_1,\dots,\om_m}$ are linearly independent forms and
$\b_{ij}\wedge\om_j=0$, Cartan's Lemma implies that
$\b_{ij}=h_{ijk}\om_k$ where $h_{ijk}=h_{ikj}$. Since $\b_{ij}$ is
symmetric, $h_{ijk}=h_{jik}$ also holds and so $h_{ijk}$  are
fully symmetric functions on the bundle $P_L$.

\vskip.5cm \noindent{\bf c. The fundamental cubic}

Let $L\subset M=\C^m$ be a SL submanifold and let $\nu\to L$ be
the normal bundle of $L$ in $M$, such that $TM\mid_L=TL\oplus\nu$.
The second fundamental form of $L$ is a quadratic form with values
in the normal bundle $\nu$ and it can be interpreted as a
traceless symmetric cubic form in the following way. The second
fundamental form $B$ of $L$ in $M$ can be written as $
B=Je_i\otimes h_{ijk}\om_j\om_k$, where $h_{ijk}$ are the fully
symmetric functions determined by $\b_{ij}$ as described above.
All the information of the second fundamental form is contained in
the symmetric cubic form $C=h_{ijk}\om_i\om_j\om_k $ which is
called the {\it fundamental cubic} of the special Lagrangian
submanifold $L$. We note that this cubic is traceless with respect
to the induced metric on $L$, $g=\om_1^2+\dots+\om_m^2$, since:
$$ tr_g C=h_{iik}\om_k=\b_{ii}=0.$$

The following result tells us that the necessary and sufficient
conditions for the existence of a special Lagrangian in $\C^m$
with a given metric and a given fundamental cubic are the Gauss
and Codazzi-type equations (\ref{Gauss}) and (\ref{Codazzi}).

Let $(L,g)\subset \C^m$ be a simply connected Riemannian manifold
of dimension $m$ and $C$ a symmetric cubic which is traceless with
respect to $g$. Choose a $g$-orthonormal coframing $\om=(\om_i)$
on an open neighborhood $U\subset L$ and define $\eta_i=0$. Now,
let $\al_{ij}=-\al_{ji}$ be the unique 1-forms on $U$ s.t.
$d\om_i=-\al_{ij}\wedge\om_j $. Write the cubic as
$C=h_{ijk}\om_i\om_j\om_k $ and set $\b_{ij}=h_{ijk}\om_k$ .

\begin{theorem}
\label{GCth}
 {\bf (see $\cite {br1}$)} Suppose that the forms $\b_{ij}$ determined by $C$ together with the
forms $\al_{ij}$ determined by $(\om_i)$ satisfy the Gauss and
Codazzi equations $(\ref{Gauss})$, $(\ref{Codazzi})$. Then there
is an isometric immersion of $(L,g)$ into $\C^m$ as a special
Lagrangian submanifold inducing $C$ as its fundamental cubic.
Moreover, this isometric immersion is {\it unique} up to rigid
motions.
\end{theorem}

%%%%%%%%%%%%%%%%%%%%%%%%%%%%%%%%%%%%%%%%%%
%%%%%%Discrete subgroups of SO(4)
%%%%%%%%%%%%%%%%%%%%%%%%%%%%%%%%%%%%%%%%%
\section{Second Order Families }
\subsection{Discrete subgroups of $\so $}

As we have seen in section 2.2.c, the second fundamental form of a
special Lagrangian submanifold $L\subset\C^m$ can be regarded as a
symmetric cubic form with vanishing trace with respect to the
induced metric $g$. It is easy to see that the symmetric cubic is
traceless if and only if it is a {\it harmonic} cubic. Therefore,
the fundamental cubic of a special Lagrangian submanifold in
$\C^m$ belongs to the space ${\mathcal H}_3(\R^4)$ of harmonic
cubics in 4 variables. This space is an irreducible $\so$-module
of dimension 16 and the action is given by
\begin{equation}
\label{action}
 (A\cdot P)x=P(xA)
\end{equation}
  where $$A=(a_{ij})\in
\so, \quad P(x)\in {\mathcal H}_3(\R^4),\quad
x=(x_1,x_2,x_3,x_4)\in\R^4$$ and
$$(xA)_i=x_j a_{ji}$$ is given by usual matrix multiplication.

We want to study the families of special Lagrangian 4-folds in
$\C^4$ whose fundamental cubic at a generic point has nontrivial
$\so$-stabilizer. The stabilizer $G$ of a polynomial $P(x)\in
{\mathcal H}_3(\R^4)$ is defined as
$$G=\{A\in\so|\ (A\cdot P)(x)=P(x), \ \mbox{for any
} x\in\R^4\}$$

 The stabilizer can be either a positive
dimensional subgroup of $\so$ or else a discrete subgroup of
$\so$. In our analysis, we need to know which are the discrete
subgroups of $\so$ that can stabilize a harmonic cubic in 4
variables.

We start by listing the discrete subgroups of $ \so$ not
containing the central rotation $-I_4$, since a subgroup of $\so$
that contains $ -I_4 $ cannot stabilize any nontrivial cubic
polynomial. For a complete proof of the classification of the
discrete subgroups of $\so$ the reader might want to consult
\cite{pv}.

In the study of the discrete subgroups of $\so$, we are going to
use the quaternionic field $\H$. Let $\E_4$ be the Euclidean
4-dimensional space and let $\{{\bf 1},{\bf i},{\bf j},{\bf k\}}$
be an orthonormal basis. We define a multiplication of elements of
$ \E_4 $ by the well-known rules: $ {\bf i}^2={\bf j}^2={\bf
k}^2=-{\bf 1},\ {\bf i}{\bf j}={\bf k},\ {\bf j}{\bf k}={\bf i},\
{\bf k}{\bf i}={\bf j} $. The elements of $\E_4$ form the
non-commutative field of quaternions $\H$. We will denote a
quaternion by the ordered set $(w,x,y,z)$ or by $w+x{\bf i}+y{\bf
j}+z{\bf k}$. For a quaternion $q=w+x{\bf i}+y{\bf j}+z{\bf k}$
 we define the conjugate $\bar q=w-x{\bf i}-y{\bf j}-z{\bf k}$ and the modulus of $q$ as
$|q|=(q \bar q)^{\frac{1}{2}}$. If $|q|=1$, we call $q$ a unit
quaternion and $U=\{q\in \H\mid |q|=1\}$ is a multiplicative group
called the group of unit quaternions.

 For any $q\in \E_4$, we define the {\it right multiplication map}
 $\ro_q \colon \E_4 \to \E_4$ by $\ro_q(x)=x \bar q $ and the {\it left
 multiplication map} $ \la_q \colon \E_4 \to \E_4$ by
 $\la_q(x)=qx$.
If $u\in U$, both $\ro_u$ and $\la_u$ are seen to
 be  in $\so$ and
 they are called the right rotation and the left rotation,
 respectively. The right rotations $\{\ro_u:u\in U\}$ form a group
 $U_+$ and the left rotations $\{\la_u:u\in U\}$ form a group $U_-$
 and both $U_+$ and $U_-$ are subgroups of $\so$, isomorphic to the
 unit quaternions group $ U$. These are different subgroups of $\so$
 and we notice that $U_+ \cap U_- = \{\pm 1 \}$.

 Consider now the homomorphism $\Phi\colon U\times U\to
 \so$ with $$\Phi(u_1,u_2)(x)=u_1x \bar
 u_2=\la_{u_1}\ro_{u_2}(x).$$ It is well-known that $\Phi$ is surjective and its kernel is the 2-element
 group $\{(1,1),(-1,-1)\} \cong \Z_2$. So, $U \times
 U/_{\{(1,1),(-1,-1)\}}\cong\so$ and to study the subgroups of $ \so $ we
 would be interested in the subgroups of the unit quaternions group
 $U$. We also define the surjective 2:1 homomorphism $\Psi\colon U\to SO(3)$ by
 $\Psi(u)(x)=ux\bar u $, whose kernel is $\{\pm 1\}$.

 It is well-known \cite{nst} that the discrete subgroups of $SO(3)$ are
 the cyclic groups, the dihedral groups and the pure symmetry groups
 of the platonic solids which are the tetrahedral, the octahedral and
 the icosahedral groups. The tetrahedral group ${\mathcal T}$ is the group of rotational
 symmetries of a tetrahedron, it has order 12 and it is isomorphic to
 $A_4$, the group of even permutations of 4 elements. The octahedral group ${\mathcal O}$
 is the group of rotational symmetries
 of an octahedral (or a cube) and it is a group of order 24, isomorphic
 to $S_4$, the permutation group of 4 elements. The icosahedral group ${\mathcal I}$
 is the group of rotational
 symmetries of a icosahedral (or a dodecahedral) and it has order 60,
 isomorphic to $A_5$.

Using the homomorphism $\Psi$, it is not hard to determine the
discrete subgroups of the group of unit quaternions $U$. A
complete description of how this is done can be found in
\cite{pv}.

\begin {proposition}
 Every finite subgroup of the unit quaternionic group $U$ is
 conjugate to one of the following groups:
\vskip-.6cm
 \begin{align*}
 {\bf C_n}&=\l\{\cos\frac{2k\pi}{n}+{\bf
 k}\sin\frac{2k\pi}{n},k=0...n-1\r\}\\
 {\bf D_n}&={\bf C_{2n}}\cup {\bf i}{\bf C_{2n}}\\
{\bf
 T}&=\ma\cup_{k=0}^{2}\l(\frac{1}{2}+\frac{1}{2}{\bf i}+\frac{1}{2}{\bf
 j}+\frac{1}{2}{\bf k}\r)^k{\bf D_2}\\
{\bf O}&={\bf T}\cup\frac{1}{\sqrt2}({\bf 1}+{\bf i}){\bf T}\\
{\bf I}&=\ma\cup_{k=0}^{4}\l(\frac{1}{2\tau}+\frac{\tau}{2}{\bf
i}+\frac{1}{2}{\bf j}\r)^k{\bf T}
\end{align*} where $\tau=\frac{\sqrt5+1}{2}$.
\end{proposition}
The binary polyhedral subgroups ${\bf T},{\bf O}$ and ${\bf I}$
are called, respectively, the binary tetrahedral, octahedral and
icosahedral subgroups of $U$ and are twice the order of the
corresponding polyhedral subgroup of $SO(3)$.

 Now, to find the discrete subgroups of $\so$ we use the 2:1
 homomorphism:
 \begin{align*}
\Phi&\colon U\ti U \to \so\\
\Phi&(l,r)(x)=l x \bar r
\end{align*}
with Ker$\,\Phi=\{(1,1),(-1,-1)\}$. We denote by
$[l,r]=\{(l,r),(-l,-r)\}$. Then we have an isomorphism $\phi\colon
U \ti U/_{\Z_2}\cong\so \mbox{ with }\phi([l,r])=\{x\to l x \bar
r\}.$

 If $\sigma\subset\so$ is a discrete subgroup, then $$L=\{l\in
 U | x\to l x \bar r \in \sigma\} \mbox{  and  } R=\{r\in
 U | x\to l x \bar r \in \sigma\}$$ are subgroups of $U$. We notice that
 $\sigma\subseteq \Phi(L\ti R)$ but the equality might not hold since it
 might be possible to find a  pair $(l,r)\in L\ti R$ such that $ x\to
 l x \bar r \not\in \sigma$. We define the subgroups
\begin{align*}
 L'=&\{l\in
 L | x\to l x\in \sigma\}=\{l\in L | (l,1)\in\sigma\}\\
 R'=&\{r\in
 R | x\to  x \bar r \in \sigma\}=\{r\in R | (1,r) \in\sigma\}
\end{align*}
  and
 there is an isomorphism between the quotient groups $L/_{L'}$ and
 $R/_{R'}$ given by $\psi(lL')=rR'$ such that
 $(l,r)\in\sigma$. The subgroup $\Phi(L'\ti R')$ is normal in $\sigma$
 and $\sigma/_{\Phi(L'\ti R')}\cong L/_{L'}$ and $R/_{R'}$. The
 subgroup $\sigma$ depends on the isomorphism $\psi$ between the
 quotient groups $L/_{L'}$ and $R/_{R'}$, different isomorphism
 possibly yielding non-conjugate subgroups in $\so$. We will denote
 the group $\sigma$ by $(L/_{L'};R/_{R'})_\psi$.

 For example, ${\bf C_m}$ is a normal subgroup of ${\bf C_{mr}}$ and
 the quotient group ${\bf C_{mr}}/_{\bf C_m}\cong \Z_r$. The elements
 of ${\bf C_{mr}}/_{\bf C_m}$ are the cosets ${\bf p}^i{\bf C_m}, i=0..r-1$,
 where ${\bf p}=\cos\frac{2\pi}{mr}+{\bf k}\sin\frac{2\pi}{mr}$ is a
 generator of ${\bf C_{mr}}$. If ${\bf q}=\cos\frac{2\pi}{nr}+{\bf
 k}\sin\frac{2\pi}{nr}$ is a generator of $ {\bf C_{nr}}$, we have the
 isomorphism $\psi_s\colon {\bf C_{mr}}/_{\bf C_m}\to {\bf
 C_{nr}}/_{\bf C_n}$ defined by $\psi_s({\bf p}^i{\bf C_m})={\bf
 q}^{si}{\bf C_n},i=0..r-1$. For each $s$ such that $(s,r)=1$ and
 $s<\frac{1}{2}r$, we get an isomorphism $\psi_s$ that gives distinct
 subgroups $({\bf C_{mr}}/_{\bf C_m};{\bf C_{nr}}/_{\bf
 C_n})_{\psi_s}. $ The subgroups of $\so$ of this form that do not
 contain the central rotation are seen to be $ ({\bf C_{2mr}}/_{\bf C_m};{\bf C_{2nr}}/_{\bf
 C_n})_{\psi_s}$, of order $mnr$ with $m$ and $n$ odd. Also, extending the
 isomorphism between ${\bf C_{2mr}}/_{\bf C_m}$ and ${\bf C_{2nr}}/_{\bf
 C_n}$ to one between ${\bf D_{mr}}/{\bf C_m}=({\bf C_{2mr}}\oplus {\bf i}{\bf
 C_{2mr}})/{\bf C_m}$ and ${\bf D_{nr}}/{\bf C_n}=({\bf C_{2nr}}\oplus {\bf i}{\bf
 C_{2nr}})/{\bf C_n}$ by $\psi_s({\bf i}{\bf p}^j{\bf C_m})={\bf
 i}{\bf q}^{sj}{\bf C_n}, j=0...r-1$, we obtain the subgroup $({\bf
 D_{mr}}/{\bf C_m};{\bf
 D_{nr}}/{\bf C_n})_{\psi_s}$ of order $2mnr$, where $m$ and $n$ are odd.

 Another subgroup of $\so$, of order 12, that does not contain the central rotation
 is $$\T=({\bf T/C_1};{\bf T/C_1})=({\bf T};{\bf T})=\{[t,t] \mid t \in{\bf
 T}\}$$

 In the case case when $L=R={\bf O}$ and $L'=R'={\bf C_1}$, we obtain
 two non-conjugate groups, depending on the automorphism of ${\bf O}$
 considered. If we take $\psi \colon {\bf O} \to{\bf O}$ to be the
 identical automorphism we obtain the subgroup $\O=({\bf O};{\bf
 O})=\{[o,o],o\in{\bf O}\}$ of order 24. If we consider the
 automorphism $\psi\colon {\bf O}={\bf T}\oplus({\bf 1}+\frac{1}{\sqrt2}{\bf i}){\bf
 T}\to {\bf O}$ with $\psi(o)=o$, if $o\in {\bf T}$ and $\psi(o)=-o$, if
 $o\in \frac{1}{\sqrt2}({\bf 1}+{\bf i}){\bf
 T}$, we obtain a different subgroup $\O^+=\{[o,o],o\in{\bf T}
 \mbox{ and }[o,-o], o\in \frac{1}{\sqrt2}({\bf 1}+{\bf i}){\bf
 T}\}$, of order 24.

 In the case when $L=R={\bf I}$ and $L'=R'={\bf C_1}$, we obtain again
 two non-conjugate subgroups of $\so$ which do not contain the central
 rotation. If we take $\psi \colon {\bf I} \to{\bf I}$ to be the
 identical automorphism, we obtain the subgroup $\I=({\bf I};{\bf
 I})=\{[l,l],l\in{\bf I}\}$, of order 60. But we notice that all the
 elements of ${\bf I}$ are in the field of rational numbers over
 $\sqrt5$ and the change of sign of $\sqrt5$ interchanges $\pm \tau$
 with $\mp \tau^{-1}$. If ${\bf p}\in {\bf I}$ is a quaternion we
 denote by ${\bf p^+}$ its image under this automorphism. Then ${\bf
 I}$ is interchanged with a group ${\bf I^+}=\ma\cup_{k=0}^{4}(\frac{\tau}{2}+\frac{1}{2\tau}{\bf
i}+\frac{1}{2}{\bf j})^k{\bf T}$ and the two groups have in common
 ${\bf T}$. If we now consider the isomorphism $\psi \colon{\bf
 I^+}\to{\bf I}, \psi({\bf p^+})={\bf p}$ then we obtain a different
 subgroup $\I^+=\{[{\bf r^+},{\bf r}],{\bf r}\in{\bf I}\}$, of order
 60. This group leaves no axis fixed and it can be shown to be the
 rotational symmetry group of a regular simplex in ${\E_4}$, with
 vertices at ${\bf 1 }$ and $\frac{1}{4}(-{\bf 1}\pm \sqrt5{\bf i}\pm \sqrt5{\bf j}
 \pm \sqrt5{\bf k}).$

 To conclude, we have the following:
\begin{proposition}
 The discrete subgroups of $\so$ that do not contain the central
 symmetry $-I_4$ are the following:
\begin{align*}
 & {\bf 1. } \ ({\bf C_{2mr}}/_{\bf C_m};{\bf C_{2nr}}/_{\bf
 C_n})_{\psi_s} \mbox{, of order  $ mnr$, where $m$, $n$  odd;}\\
 &{\bf 2. } \ ({\bf
  D_{mr}}/{\bf C_m};{\bf
 D_{nr}}/{\bf C_n})_{\psi_s} \mbox{, of order  $2mnr$, where $m$, $n$
 odd;}\\
 &{\bf 3. } \ \T=({\bf T}/{\bf C_1};{\bf T}/{\bf C_1})=\{[t,t]\mid t\in{\bf T}\}
 \mbox{, of order $12$;}\\
&{\bf 4. } \ \O=({\bf O}/{\bf C_1};{\bf O}/{\bf C_1})=\{[o,o]\mid
o\in{\bf O}\}
 \mbox{, of order $24$;}\\
&{\bf 5. } \ \O^+=({\bf O}/{\bf C_1};{\bf O }/{\bf
C_1})=\{[o,o]\mid o\in{\bf T} \mbox{
  and } [o,-o],o\in \frac{1}{\sqrt2}({\bf 1}+{\bf i}){\bf T}\}
 \mbox{, of order $24$;}\\
&{\bf 6. } \ \I=({\bf I}/{\bf C_1};{\bf I}/{\bf C_1})=\{[l,l]\mid
l\in{\bf I}\}
 \mbox{, of order $60$;}\\
&{\bf 7. } \ \I^+=({\bf I^+}/{\bf C_1};{\bf I}|{\bf
C_1})=\{[r^+,r]\mid r\in{\bf I}\}
 \mbox{, of order $60$.}
\end{align*}
\end{proposition}
 For complete proof of this proposition, the reader should consult
 \cite {pv} and {\cite {meg}.

%%%%%%%%%%%%%%%%%%%%%%%%%%%%%%%%%%%%%%%%%%
%%%%%%Continuous symmetry
%%%%%%%%%%%%%%%%%%%%%%%%%%%%%%%%%%%%%%%%%
\subsection{Continuous symmetry}

Any maximal torus in $\so$ is conjugate to the group:
$$\l\{\begin{pmatrix}
e^{i\theta_1}&0\\0&e^{i\theta_2}\end{pmatrix},\
\theta_1,\theta_2\in [0,2\pi)\r\}$$ We are looking to determine
the orbits of the action (\ref{action}) that have non-trivial
stabilizer under the action of $\so$. In what follows, a positive
dimensional stabilizer will be called a continuous stabilizer. The
special Lagrangian 4-folds whose fundamental cubic at each point
has a continuous stabilizer will be said to have {\it continuous
symmetry}. If the stabilizer is a finite subgroup of $\so$, we
will say that the special Lagrangian 4-fold has {\it discrete
symmetry}.

First, we are going to classify the harmonic cubic polynomials in
4 variables $\{x_1,x_2,x_3,x_4\}$ whose stabilizer is a continuous
subgroup of $\so$.

\begin{proposition}
\label{cont}
 The $\so$-stabilizer of $h\in{\mathcal H}_3(\R^4)$ is a continuous subgroup
of $\so$ if and only if $h$ lies on the $\so$-orbit of exactly one
of the following polynomials:
\\

\non $1.\ 0\in {\mathcal H}_3(\R^4),\mbox{whose stabilizer
is}\so;$\vskip.2cm

\non $2.\ rx_1(x_1^2-x_2^2-x_3^2-x_4^2) \mbox{ for some } r>0,
\mbox{ whose stabilizer is the subgroup $ SO(3)$, consisting}$
\vskip.2cm \non $\mbox{ of rotations in the $3$-space
$(x_2,x_3,x_4)$};$\vskip.2cm

\non $3.\ r[(x_1^2-x_2^2)x_3+2x_1x_2x_4], \mbox{ for some } r>0,
\mbox{ whose stabilizer is the subgroup $O(2)$
generated}$\vskip.2cm \non $\mbox{ by rotations by an arbitrary
angle in the $(x_1,x_2)$-plane and twice that angle in the
}$\vskip.2cm \non $\mbox{ $(x_3,x_4)$-plane and the element
$\biggl(\begin{smallmatrix}1&0&0&0\\0&-1&0&0\\0&0&1&0\\0&0&0&-1\end{smallmatrix}\biggl)$};$\vskip.2cm

\non $4.\ r(x_1^3-3x_1x_2^2) \mbox { for some } r>0, \mbox{whose
 stabilizer is the subgroup $ SO(2)\ltimes S_3$, where $S_3$  is
 the}$\vskip.2cm
\non $\mbox{ symmetric group on $3$ elements generated by the
rotation by an angle of $\frac{2\pi}{3}$  in the
$(x_1,x_2)$-}$\vskip.2cm \non $\mbox{ plane and the element
$\biggl(\begin{smallmatrix}1&0&0&0\\0&-1&0&0\\0&0&-1&0\\0&0&0&1\end{smallmatrix}\biggl)$,
 and $SO(2)$ is the group of rotations in the $(x_3,x_4)$-plane;
 }$\vskip.2cm

\non $5.\  r(x_1^3-3x_1x_2^2)+3vx_1(x_1^2+x_2^2-2x_3^2-2x_4^2)
\mbox { for some $v>0$ satisfying $ r \not = 3v$  whose
stabilizer},$\vskip.2cm \non $\mbox {is the $O(2)$-subgroup
generated by rotations in
 the $(x_3,x_4)$-plane and the element
 $\biggl(\begin{smallmatrix}1&0&0&0\\0&-1&0&0\\0&0&1&0\\0&0&0&-1\end{smallmatrix}\biggl)$;
 }$\vskip.2cm

\non $6.\
r(x_1^3-3x_1x_2^2)+s(3x_1^2x_2-x_2^3)+3vx_1(x_1^2+x_2^2-2x_3^2-2x_4^2)\mbox
{ for
 some $s>0$  and $v>0$ whose}$\vskip.2cm
\non$\mbox{
 stabilizer is the $SO(2)$-subgroup generated by rotations in the
 $(x_3,x_4)$-plane.}$
\end{proposition}

 \noindent{\it Remark }(special values). The case $r=3v$ of case 5 above reduces to case 2, when the
stabilizer is $SO(3)$.

 \pf Suppose $h\in{\mathcal H}_3(\R^4)$ has a nontrivial stabilizer
$G\subseteq\so$. If $G=\so$, then $h=0$  since ${\mathcal
H}_3(\R^4)$ is an irreducible representation of $\so$. We suppose
from now on that $h\not= 0$. Being a stabilizer, $G$ is a closed
subgroup of $\so$, therefore it is compact and has a finite number
of components.

We suppose $G$ is not discrete. Then, its identity component $H$
is a closed connected subgroup of $\so$. The algebra
$\mathfrak{h}$ of $H$ is a subalgebra of ${\mathfrak so}(4)$.
Using the 2:1 homomorphism $\Phi\colon U\times U\to
 \so$, $\Phi(u_1,u_2)(x)=u_1x \bar
 u_2$ from section 3.1., it is easy to see that ${\mathfrak so}(4) \cong
{\mathfrak so}(3)_+ \oplus {\mathfrak so}(3)_-$, where $
{\mathfrak so}(3)_+ $ and ${\mathfrak so}(3)_- $ are two different
copies of ${\mathfrak so}(3)$ with intersection the 0 vector.
Since $\dim {\mathfrak so}(4)=6$, there are the following
possibilities for the subalgebra $\mathfrak{h}$:

1) dim ${\mathfrak h}=5$. This is not possible for the following
reason:
  ${\mathfrak h}\cap {\mathfrak so}(3)_+ \subset {\mathfrak so}(3)_+$ is a
  subalgebra of dimension at least 2 and ${\mathfrak so}(3)_+$ has no
  subalgebras of dimension 2. Therefore,  ${\mathfrak h}\cap {\mathfrak so}(3)_+={\mathfrak
  so}(3)_+$ which implies that ${\mathfrak h}\supseteq {\mathfrak so}(3)_+$. Similarly,
  it can be shown that ${\mathfrak h}\supseteq {\mathfrak so}(3)_-$ and it follows from here that
  ${\mathfrak h}={\mathfrak so}(4)$, which gives a contradiction.

2) dim ${\mathfrak h}=4$. Then ${\mathfrak h}_+={\mathfrak h}\cap
{\mathfrak so}(3)_+$
  is an ideal of dimension at least 1 in ${\mathfrak so}(3)_+$ and ${\mathfrak h}_-=
  {\mathfrak h}\cap {\mathfrak so}(3)_-$
  is an ideal of dimension at least 1 in ${\mathfrak so}(3)_-$. Consider the
  projections $\pi_{\pm}\colon{\mathfrak h}\to {\mathfrak
  so}(3)_{\pm}$. Since Ker $\pi_{\pm}= {\mathfrak h}_\mp$,  it follows that
  Ker $\pi_{\pm}$ can have dimension 1 or 3. If Ker $\pi_-$ has dimension 1, $\pi_-$ is onto and
  Ker $\pi_+$ has dimension 3. In this case, it follows that ${\mathfrak h}_+ \cong {\mathfrak
  so}(2)_+$ and we obtain that ${\mathfrak h}={\mathfrak so}(2)_+ \oplus
  {\mathfrak so}(3)_-$. If the dimension of Ker $\pi_-$ is 3, then the
  dimension of Ker $\pi_+$ is 1 and we obtain ${\mathfrak
  h}={\mathfrak so}(3)_+ \oplus {\mathfrak so}(2)_-$. But ${\mathfrak
  so}(3)_\pm$ acts like the group $SU(2)$ on the space of complexified harmonic cubics
  in four variables $\{z_1,z_2,{\bar z_1},{\bar z_2}\}$, where $z_1=x_1+ix_2, z_2=x_3+ix_4$ and calculations show that
  this action does not preserve any nontrivial element. Therefore, in this case
  $H$ does not preserve any nontrivial harmonic cubic.

3) dim ${\mathfrak h}=3$. In this case, one can show that, up to
conjugacy, the only
  possibilities for ${\mathfrak h}$ are ${\mathfrak so}(3)_+$,
  ${\mathfrak so}(3)_-$ and diag$({\mathfrak so}(3)_+ \oplus
  {\mathfrak so}(3)_-)$. But, as discussed in 2) above, ${\mathfrak so}(3)_\pm$
  does not preserve any cubic polynomial in 4 variables and
  consequently we can discard these cases.

 We study now the case $${\mathfrak h}=\mbox{diag}({\mathfrak so}(3)_+ \oplus
  {\mathfrak so}(3)_-)=\{x_++x_-,\ x_+\in {\mathfrak so}(3)_+,\ x_-\in
  {\mathfrak so}(3)_-\} $$ We can see that, up to conjugacy,
 $$ \mbox {diag}({\mathfrak so}(3)_+ \oplus
  {\mathfrak so}(3)_-)=\l\{
\begin{pmatrix} 0 & 0 & 0 & 0 \\ 0 & 0& -2c & 2b\\ 0  & 2c& 0 & -2a \\
0 & -2b & 2a & 0
\end{pmatrix},a,b,c\in\R \r\}$$
Therefore, $G$ can be  either one of the groups: $G=\l\{
\begin{pmatrix} 1 & 0\\ 0 & A
\end{pmatrix}, A\in SO(3) \r\}$ or\\
$G=\l\{
\begin{pmatrix} {\mbox det}(A) & 0 \\ 0 & A
\end{pmatrix}, \ A\in O(3) \r\}$.
We can easily see that the cubic polynomials fixed by the identity
component $H$ are linear combinations of $x_1^3$ and $
x_1(x_2^2+x_3^2+x_4^2)$. It is obvious that the only combination
that would make the polynomial harmonic is $
P=rx_1(x_1^2-x_2^2-x_3^2-x_4^2)$, for some $r\not=0$. One can
verify that the full stabilizer of $P$ is $SO(3)$. By a rotation
that reverses the $x_1$-axis, if necessary, we can assume that
$r>0$.

4) dim ${\mathfrak h}=2$. In this case, ${\mathfrak h}={\mathfrak
so}(2)_+ \oplus
  {\mathfrak so}(2)_-$ and $H$ is conjugate to the maximal torus $H=\l\{
\begin{pmatrix} e^{i\theta_1} & 0\\ 0 & e^{i\theta_2}
\end{pmatrix},\ \theta_1,\theta_2\in [0,2\pi)\r\}$. It is easy to
see that $H$ does not stabilize any symmetric cubic in 4
variables.

5) dim ${\mathfrak h}=1$. In this case, one can show that the only
  1-dimensional ideals in ${\mathfrak so}(4)$ are conjugate to:
$$\hf_{p,q}=\{(px,qx)|\ x\in\sof(2),\ p,q\in\Z, \ (p,q)=1\}$$
It follows that the identity component
$H_{p,q}=\l\{\begin{pmatrix} R(p\theta) & 0\\ 0 &
R(q\theta)\end{pmatrix},\theta\in\R\r\}$ consists of rotations of
angle $p\theta$ in the $(x_1,x_2)$-plane and of angle $q\theta$ in
the $(x_3,x_4)$-plane. We are looking for harmonic cubics in
$\{x_1,x_2,x_3,x_4\}$ preserved by $H_{p,q}$, for some $p$ and $q$
integers.

Let $V_n$ be the irreducible representation of $SO(2)$ given by
rotations of speed $n$: $e^{i\theta}.z=e^{ni\theta}z$, where
$z\in\C$. In our case, the speed $p$ representation $V_p$ is given
by the action on the $(x_1,x_2)$-plane:
$e^{i\theta}.z_1=e^{ip\theta}z_1$ and $V_q$ is given by the action
on the $(x_3,x_4)$-plane: $e^{i\theta}.z_2=e^{iq\theta}z_2$, where
$z_1=x_1+ix_2$ and $z_2=x_3+ix_4$.

Under the action of $H_{p,q}$, the space of symmetric polynomials
in 4 variables ${\Scal}^3(\R^4)$ decomposes as: $$
{\Scal}^3(V_p\oplus
  V_q)={\Scal}^3(V_p)\oplus ({\mathcal S}^2(V_p)\otimes {\Scal}^1(V_q))\oplus
  ({\Scal}^1(V_p)\otimes {\Scal}^2(V_q))\oplus {\Scal}^3(V_q).$$ But one
can see that ${\Scal}^3(V_p)\cong V_{3p}\oplus V_p$. A basis in
$V_{3p}$ is $\{\Re z_1^3,\Im
z_1^3\}=\{x_1^3-3x_1x_2^2,3x_1^2x_2-x_2^3\}$ and a basis in $V_p$
is $\{\Re( z_1 \bar z_1 z_1),\Im(z_1 \bar z_1
z_1)\}=\{(x_1^2+x_2^2)x_1,(x_1^2+x_2^2)x_2\}$. Similarly,
${\Scal}^2(V_p)\cong V_{2p}\oplus V_0\cong V_{2p}\oplus \R$ and we
calculate that:
\begin{align*}
{\Scal}^3(\R^4)=&(V_{3p}\oplus V_p)\oplus((V_{2p}\oplus\R)\otimes
V_q)\oplus (V_p\otimes
(V_{2q}\oplus \R)) \oplus(V_{3q}\oplus V_q)\\
=&V_{3p}\oplus V_p\oplus V_{2p+q}\oplus V_{2p-q}\oplus
V_{p+2q}\oplus V_{p-2q}\oplus V_{3q}\oplus V_q
\end{align*}
 where we used that $V_n\otimes V_m=V_{n+m}\oplus V_{n-m}$. This decomposition is
irreducible and the action has a fixed vector if and only if one
of the $V$'s in the above direct sum is $V_0$, the 1-dimensional
representation on which the action is trivial. This implies one of
the following possibilities for the values of $p$ and $q$:
$$p=0,\ q=0,\ p-2q=0,\ q-2p=0,\ p+2q=0  \mbox{ or } 2p+q=0$$ We note that $p$ and $q$ have symmetric
roles since the planes $(x_1,x_2)$ and $(x_3,x_4)$ can be
interchanged by an orthogonal transformation. Therefore, up to
conjugation with an element in $\so$, the only possible cases are:
$p=0$, $q=2p$ and $q=-2p$. In the case $q=2p$, the fact that
$(p,q)=1$ implies that $p=1$ and $q=2$ and in the case $q=-2p$, we
can take $p=1, q=-2$. The conclusion is that, unless one of these
conditions is satisfied, the group $H_{p,q}$ does not preserve any
cubic symmetric polynomial in four variables. But since the
stabilizer in $\so$ coincides with the stabilizer in $O(4)$, then
up to conjugacy in $O(4)$, the last two cases are the same. We
will study each of these cases separately.

a) $p=1$ and $q=2$. Then $H_{1,2}= \l\{\begin{pmatrix} e^{i\theta}
& 0 \\ 0 &
 e^{2i\theta}\end{pmatrix}, \theta\in [0,2\pi)\r\}$. If $P$ is a complexified
 polynomial in the variables
 $\{z_1,z_2,{\bar z_1},{\bar z_2}\}$, fixed by $H_{1,2}$, it is easy to
 see that $P$ should lie in $V_2\otimes V_1$, therefore $P=az_1^2{\bar z_2}+b{\bar z_1}^2
 z_2,a,\ b\in\C$. Now, $P$ is a real harmonic polynomial if $b={\bar a}$.
 So, any real harmonic cubic $C$ preserved by this group is of the form:
 $$C=\Re(az_1^2{\bar z_2})=r[(x_1^2-x_2^2)x_3+2x_1x_2x_4]+s[2x_1x_2x_3-(x_1^2-x_2^2)x_4],$$
 with $r,s\in \R$. If $r^2+s^2\not= 0$, by applying a rotation of angle
 $\arctan(\frac{s}{r})$ in the $(x_3,x_4)$-plane, we can assume that
 $s=0$. We can also assume that $r\geq 0 $. The conclusion is that all
 the harmonic cubics in 4 variables stabilized by $H_{1,2}$ are on the
 $\so$-orbit of the cubic $h=r[(x_1^2-x_2^2)x_3+2x_1x_2x_4]$. The full stabilizer of $h$ can be shown to be the
 disconnected 2-piece subgroup $H_{1,2}\cup gH_{1,2}\subset \so$, where
 $g=\Biggl(\begin{smallmatrix}1&0&0&0\\0&-1&0&0\\0&0&1&0\\0&0&0&-1\end{smallmatrix}\Biggl)$.
 It is isomorphic to $O(2)$, because $O(2)$ is the only non-abelian
 2-component compact group of dimension one.

b) If $p=0$, we can consider $q=1$. In this case, the identity
 component of $G$ is $H_{0,1}= \l\{\begin{pmatrix} I_2 & 0 \\ 0 &
 e^{i\theta}\end{pmatrix}, \theta\in\R\r\}\cong S^1$. A complexified cubic polynomial $C$, fixed by
 this group, should belong to $V_0\otimes V_1$, therefore it should be a
 linear combination of $z_1z_2{\bar z_2}$, $z_1^3$ and $z_1^2{\bar
 z_1}$.
Calculations show $C$ is harmonic if and only if it is a linear
combination of the harmonic polynomials $\{z_1^3, z_1^2{\bar
z_1}-2z_1z_2{\bar z_2}\}$. Therefore, the fixed real harmonic
cubic polynomials $C$ in 4 variables are of the form
$$C=r(x_1^3-3x_1x_2^2)+s(3x_1^2x_2-x_2^3)+vx_1(x_1^2+x_2^2-2x_3^2-2x_4^2)+
ux_2(x_1^2+x_2^2-2x_3^2-2x_4^2),$$ with $r,s,u,v\in\R$. By making
a rotation in the $(x_1,x_2)$-plane, we can suppose that $u=0$.

It remains now to determine the full stabilizer of the polynomial
$$r(x_1^3-3x_1x_2^2)+s(3x_1^2x_2-x_2^3)+vx_1(x_1^2+x_2^2-2x_3^2-2x_4^2),$$
which we denote by $G$.

If $s\not = 0$ and $v\not = 0$, calculations show that the full
stabilizer of $h$ is just the identity component, so $G=SO(2)$. By
making some rotations, if necessary, we may assume that $s,\ v>0$.

If $s=0$, $v\not =0$, and $r\not = 3v$, then
$h=r(x_1^3-3x_1x_2^2)+vx_1(x_1^2+x_2^2-2x_3^2-2x_4^2)$ and the
full stabilizer is isomorphic to an $O(2)$-subgroup, because we
are also allowed to flip the signs of $x_2$ and $x_4$. In this
case we can suppose that $v>0$. In the case $s=0$ and $r = 3v$,
$h$ becomes $6vx_1(x_1^2-x_2^2-x_3^2-x_4^2)$ and we saw that this
polynomial has stabilizer $SO(3)$.

Finally, if $s=v=0$ and $r>0$, the polynomial
$h=r(x_1^3-3x_1x_2^2)$ is preserved by the identity component
$S^1$, but we can see that it is also fixed by the element of
order 3, $A=\bigl(\begin{smallmatrix}e^{\frac{2\pi i}{3}}& 0\\ 0 &
I_2
\end{smallmatrix}\bigl)$ and by the element of order 2, $g=\Biggl(\begin{smallmatrix}
1&0&0&0\\0&-1&0&0\\0&0&-1&0\\0&0&0&1\end{smallmatrix}\Biggl)$ and
both $A$ and $g$ do not belong to the identity component. The
rotation $A$ and the reflection $g$ form a group isomorphic to
$S_3$, the symmetric group in 3 elements. $S_3$ acts on the
identity component $S^1$ by conjugation and it is easy to compute
that the action of $A$ on $S^1$ is trivial, while $g$ acts by
flipping the circle $S^1$, namely $g.\begin{pmatrix} I_2 & 0 \\ 0
&
 e^{i\theta}\end{pmatrix}=\begin{pmatrix} I_2 & 0 \\ 0 &
 e^{-i\theta}\end{pmatrix}$. Using the exact sequence
$$
\begin{array}{ccccccccc}
0 &\ra &S_1 &\ra & G & \ra & S_3 &\ra & 0
\end{array}
$$
one can verify that $G=S^1\ltimes S_3$. This completes the proof
of Prop \ref{cont}. $\Box$

\vskip.2cm

In the next step, we are going to analyze each of the cases given
by Prop \ref{cont} and classify the SL 4-folds in $\C^4$ whose
fundamental cubic form has a continuous symmetry at each point.

\subsubsection{SO(3)-symmetry}

{\it Example 1.} In their paper \cite{hl}, Harvey and Lawson found
the following special Lagrangian submanifolds of $\C^4$, invariant
under the diagonal action of $SO(4)$ on $\C^4=\R^4\ti \R^4$:
\begin{equation}
\label{eqHL}
 L_c=\{(s+it){\bf u}| \ {\bf u}\in S^3\subset\R^4, \Im
(s+it)^4=c\},
\end{equation}
where $c\in \R$ is a constant. The variety $L_0$ is an
 union of four special Lagrangian 4-planes and when $c\not = 0$, each component of $L_c$ is diffeomorphic to
$\R\ti S^3$ and it is asymptotic to one pair of 4-planes in $L_0$.

\begin{theorem}

 If $L\subset \C^4$ is a connected nontrivial special Lagrangian submanifold
 whose cubic fundamental form has an $SO(3)$-symmetry at each point,
 then $L$ is, up to rigid motion, an
 open subset of one of the Harvey-Lawson examples.

\end{theorem}

\pf In the above hypotheses, a trivial special Lagrangian
submanifold is a special Lagrangian 4-plane. We can see that the
fundamental cubic $C$ of $L$ lies on the orbit of the 0 cubic if
and only $L$ is trivial. Therefore, assume the fundamental cubic
is not identically vanishing. The locus where $C$ vanishes is a
proper real-analytic subset of $L$, so its complement $L^*$ is
open and dense in $L$. By replacing $L$ by its component $L^*$, we
can assume without loss of generality that $C$ is nowhere
vanishing on $L$. Since the stabilizer of $C_x$ is $SO(3)$ for all
$x\in L$, Proposition \ref{cont} implies the existence of a
positive real-analytic function $r\colon L \to \R^+ $ with the
property that the equation
$$ C=3r\om_1(\om_1^2-\om_2^2-\om_3^2-\om_4^2)$$
defines an $SO(3)$-subbundle $F$ of the bundle $P_L$ of
$L$-adapted coframes. On the subbundle $F$, the following
identities hold:
\begin{equation}
\label{eq1}
\begin{pmatrix}\b_{11}&\b_{12}&\b_{13}&\b_{14}\\
\b_{21}&\b_{22}&\b_{23}&\b_{24}\\
\b_{31}&\b_{32}&\b_{33}&\b_{34}\\
\b_{41}&\b_{42}&\b_{43}&\b_{44}
\end{pmatrix}=\begin{pmatrix}\ 3r\om_1&-r\om_2&-r\om_3&-r\om_4\\
-r\om_2&-r\om_1&0&0\\
-r\om_3&0&-r\om_1&0\\
-r\om_4&0&0&-r\om_1
\end{pmatrix}
\end{equation}
Because $F$ is an $SO(3)$-bundle, the forms $\a_{21},a_{31}$ and
$\a_{41}$ vanish mod $\{\om_1,\om_2,\om_3,\om_4\}$, meaning that
there are functions $t_{ij}$ on $F$ such that:
\begin{align}
&\a_{21}=t_{21}\om_1+t_{22}\om_2+t_{23}\om_3+t_{24}\om_4 \notag \\
& \a_{31}=t_{31}\om_1+t_{32}\om_2+t_{33}\om_3+t_{34}\om_4 \label{eq2}  \\
&\a_{41}=t_{41}\om_1+t_{42}\om_2+t_{43}\om_3+t_{44}\om_4 \notag
\end{align}

Also, there exist functions $r_i$, $i=1,2,3,4$ on $F$ such that:
\begin{equation}
\label{eq3} dr=\sum_{i=1}^{4}r_i\om_i.
\end{equation}

Substituting the relations (\ref{eq1}), (\ref{eq2}) and
(\ref{eq3}) into the identities
\begin{equation}
\label{eq4} d\b_{ij}=-\b_{ik}\wedge\al_{kj}-\al_{ik}\wedge\b_{kj}
\end{equation}
and using the identities $d\om_i=-\a_{ij}\wedge \om_j$, one gets
polynomial relations among $r_i$, $t_{ij}$ which can be solved,
leading to relations of the form:
\begin{equation}
\label{eq5} \a_{21}=t\om_2, \quad \a_{31}=t\om_3, \quad
\a_{41}=t\om_4, \quad dr=-5rt\om_1
\end{equation}
where we denoted $t_{22}=t_{33}=t_{44}$ by $t$.

Differentiating the last equation in (\ref{eq5}), we get
$0=d(dr)=-5rd(t)\wedge\om_1$, implying that there exits a function
$u$ on $F$ such that
\begin{align}
\label{eq6} dt=u\om_1.
\end{align}
Substituting the relations (\ref{eq5}) and (\ref{eq6}) into the
identities
\begin {equation}
\label{eq7} d\a_{ij}=-\a_{ik}\wedge\a_{kj}+\b_{ik}\wedge\b_{kj}
\end{equation}
and expanding out using the identities $d\om_i=-\a_{ij}\wedge
\om_j$ implies the relations:
\begin{align*}
u&=4r^2-t^2\\
d\a_{32}&=\a_{43}\wedge\a_{42}+(t^2+r^2)\om_3\wedge\om_2\\
d\a_{42}&=\-a_{43}\wedge\a_{32}+(t^2+r^2)\om_4\wedge\om_2\\
d\a_{43}&=\a_{42}\wedge\a_{32}+(t^2+r^2)\om_4\wedge\om_3
\end{align*}
Differentiating the last equations yields only identities.

So, $F\ra L$ is a $SO(3)$-bundle on which the 1-forms
$\{\om_1,\om_2,\om_3,\om_4,\a_{32},\a_{42},\a_{43}\}$ form a basis
and they satisfy the structure equations:
\begin{align}
d\om_1&=0 \notag \\
d\om_2&=t\om_1\wedge\om_2+\a_{32}\wedge\om_3+\a_{42}\wedge\om_4
\notag
\\
d\om_3&=t\om_1\wedge\om_3-\a_{32}\wedge\om_2+\a_{43}\wedge\om_4
\notag
\\
d\om_4&=t\om_1\wedge\om_4-\a_{42}\wedge\om_2-\a_{43}\wedge\om_3
\notag\\
d\a_{32}&=\a_{43}\wedge\a_{42}+(t^2+r^2)\om_3\wedge\om_2 \label{eq8}\\
d\a_{42}&=-\a_{43}\wedge\a_{32}+(t^2+r^2)\om_4\wedge\om_2 \notag \\
d\a_{43}&=\a_{42}\wedge\a_{32}+(t^2+r^2)\om_4\wedge\om_3 \notag \\
dr&=-5rt\om_1 \notag \\
dt&=(4r^2-t^2)\om_1 \notag
\end{align}
and the exterior derivatives of these equations are identities.

The last two equations in (\ref{eq8}) imply that
$$\frac{dr}{-5rt}=\frac{dt}{4r^2-t^2}=\om_1.$$ This yields
$d(r^{\frac{8}{5}}+t^2r^{-\frac{2}{5}})=0$ and since $L$ and $F$
are connected, there exists a function $\theta$ on $L$ with
$|\theta|<{\pi\over 8}$ such that:
\begin{align}
r^{\frac{4}{5}}&=c^{\frac{4}{5}}\cos4\theta \notag \\
r^{-\frac{1}{5}}t&=c^{\frac{4}{5}}\sin4\theta \notag
\end{align}
From these last two equations and from last equation in
(\ref{eq8}), it follows that
$$\om_1=\frac{d\theta}{c (\cos4\theta)^{\frac{5}{4}}}$$
Now, setting $\eta_i=c(\cos 4\theta)^{1/4}\om_i$ for $i=2,3$ and 4
yields:
\begin{align*}
d\eta_2&=-\a_{23}\wedge\eta_3-\a_{24}\wedge\eta_4 \\
d\eta_3&=-\a_{32}\wedge\eta_2-\a_{34}\wedge\eta_4 \\
d\eta_4&=-\a_{42}\wedge\eta_2-\a_{43}\wedge\eta_3\\
d\a_{32}&=-\a_{34}\wedge\a_{42}+\eta_3\wedge\eta_2\\
d\a_{42}&=-\a_{43}\wedge\a_{32}+\eta_4\wedge\eta_2\\
d\a_{43}&=-\a_{42}\wedge\a_{23}+\eta_4\wedge\eta_3
\end{align*}
The above equations represent the structure equations of the
metric of constant curvature 1 on the 3-sphere $S^3$.

Conversely, if $d\sigma^2$ is the metric of constant curvature 1
on $S^3$, then, on the product $L=(-{\pi\over 8},{\pi\over 8})\ti
S^3$, consider the quadratic form
$$g=\frac{d\theta^2+\cos^2 {4\theta} d\sigma^2}{c^2 (\cos
{4\theta})^{5/2}}$$ and the cubic form
$$C=3\frac{\cos^2 4\theta d\sigma^2 d\theta-d\theta^3}{c^2 (\cos{4\theta})^{5/2}}$$

The pair $(g,C)$ satisfies the Gauss and Codazzi equations and by
Theorem \ref{GCth} this implies that $(L,g)$ can be isometrically
immersed as a special Lagrangian 4-fold in $\C^4$ inducing $C$ as
its fundamental cubic. For each value of $c$, there exists a
unique corresponding special Lagrangian 4-fold. Since the
structure equations (\ref{eq8}) have an $\so$-symmetry and Harvey
and Lawson \cite{hl} found that all the special Lagrangian 4-folds
in $\C^4$, invariant under the diagonal action of $\so$, can be
written explicitly as (\ref{eqHL}), the conclusion of the theorem
follows. $\Box$

\subsubsection{O(2)-symmetry}

According to Prop.\ref{cont}, there are two cases of
$O(2)$-symmetry. The first one gives the following:
\begin{theorem}
There is no connected nontrivial special Lagrangian $4$-fold in
$\C^4$ whose cubic fundamental form has an $O(2)$-symmetry at each
point, where $O(2)$ is the subgroup $ S^1\cup g_0S^1$ with
$S^1=\biggl\{ \begin{pmatrix}e^{i\theta}&0
\\0&e^{2i\theta}\end{pmatrix}, \theta\in\R \biggl\}$ and
$g_0=\l(\begin{smallmatrix}1&0&0&0\\0&-1&0&0\\0&0&1&0\\0&0&0&-1\end{smallmatrix}\r).$
\end{theorem}

\pf Let $L$ be a special Lagrangian 4-fold that satisfies the
hypotheses of the theorem and let $C$ be its fundamental cubic.
Proposition \ref{cont} implies that there exists a function
$r\colon L\to \R_+$ for which the equation
$$ C=3r[(\om_1^2-\om_2^2)\om_3+2\om_1\om_2\om_4]$$
defines an $O(2)$-subbundle $F\subset P_L$ of the $L$-adapted
coframe bundle $P_L \ra L$. On the subbundle $F$, the following
identities hold:
\begin{equation}
\label{eq10}
\begin{pmatrix}\b_{11}&\b_{12}&\b_{13}&\b_{14}\\
\b_{21}&\b_{22}&\b_{23}&\b_{24}\\
\b_{31}&\b_{32}&\b_{33}&\b_{34}\\
\b_{41}&\b_{42}&\b_{43}&\b_{44}
\end{pmatrix}=\begin{pmatrix}r\om_3&r\om_4&r\om_1&r\om_2\\
r\om_4&-r\om_3&-r\om_2&r\om_1\\
r\om_1&-r\om_2&0&0\\
r\om_2&r\om_1&0&0
\end{pmatrix}
\end{equation}

Since $F$ is an $O(2)$-bundle, the following relations hold:
$\a_{31}=\a_{41}=\a_{32}=\a_{42}=\a_{43}-2\a_{21}\equiv 0 \mod
\{\om_1,\om_2,\om_3,\om_4\}$, meaning that there exist functions
$t_{ij}$ on $F$ such that:
\begin{align}
\a_{42}&=t_{11}\om_1+t_{12}\om_2+t_{13}\om_3+t_{14}\om_4 \notag \\
\a_{32}&=t_{21}\om_1+t_{22}\om_2+t_{23}\om_3+t_{24}\om_4\notag \\
\a_{31}&=t_{31}\om_1+t_{32}\om_2+t_{33}\om_3+t_{34}\om_4 \label{eq11}\\
\a_{41}&=t_{41}\om_1+t_{42}\om_2+t_{43}\om_3+t_{44}\om_4 \notag \\
\a_{43}-2\a_{21}&=t_{51}\om_1+t_{52}\om_2+t_{53}\om_3+t_{54}\om_4
\notag
\end{align}
Moreover, there exist functions $r_i$ on $F$, $i=1,2,3$ and 4 such
that
\begin{equation}
\label{eq12} dr=\sum_{i=1}^{4}r_i\om_i.
\end{equation}
Substituting the relations (\ref{eq10}), (\ref{eq11}) and
(\ref{eq12}) into the identities
\begin{equation}
d\b_{ij}=-\b_{ik}\wedge\al_{kj}-\al_{ik}\wedge\b_{kj}
\end{equation}
and using the identities $d\om_i=-\a_{ij}\wedge \om_j$, one gets
polynomial relations among $r_i$, $t_{ij}$ which can be solved,
leading to relations of the form:
\begin{equation}
\label{eq13} \a_{31}=\a_{41}=\a_{32}=\a_{42}=\a_{43}-2\a_{21}=0,
\quad dr=0
\end{equation}
Substituting (\ref{eq10}) and (\ref{eq13}) into the identities
$d\a_{ij}=-\a_{ik}\wedge\a_{kj}+\b_{ik}\wedge\b_{kj}$ yields
$r=0$, contrary to the hypothesis. $\Box$
\smallskip

The second case of symmetry $O(2)$ yields the following partial
result:
\begin{proposition}
 There is a $2$-parameter family of connected special
Lagrangian $4$-folds with the property that the fundamental cubic
at each point is isomorphic to $O(2)$, where $O(2)$ is the
subgroup $ S^1\cup g_0S^1$ with $S^1=\l\{
\begin{pmatrix}I_2&0
\\0&e^{i\theta}\end{pmatrix}, \theta\in\R\r\}$ and
$g_0=\l(\begin{smallmatrix}1&0&0&0\\0&-1&0&0\\0&0&1&0\\0&0&0&-1\end{smallmatrix}\r)$.
\end{proposition}
\pf Let $L$ be a special Lagrangian $4$-fold that satisfies the
hypotheses of the theorem and let $C$ be its fundamental cubic.
Proposition \ref{cont} implies that there exist functions
$r,v\colon L\to \R$, $v\geq 0, r\not = 3v$ and an $O(2)$-subbundle
$F\subset P_L$ over $L$ on which the following identity holds:
\begin {equation}
C=r(\om_1^3-3\om_1\om_2^2)+3v\om_1(\om_1^2+\om_2^2-2\om_3^2-2\om_4^2).
\end{equation}

Since $F$ is an $O(2)$-bundle, the following relations hold:
$\a_{21}=\a_{31}=\a_{42}=\a_{32}=\a_{42}\equiv 0 \mod
\{\om_1,\om_2,\om_3,\om_4\}$ and doing the differential analysis
as in the previous cases we obtain the following structure
equations on the subbundle $F$:
\begin{align}
d\om_1&=(r-3v)t_2\om_1\wedge \om_2 \notag\\
d\om_2&=(r-v)t_1\om_1\wedge \om_2 \notag\\
d\om_3&=2vt_1\om_1\wedge \om_3+2vt_2 \om_2\wedge
\om_3+\a_{43}\wedge
\om_4 \notag\\
d\om_4&=2vt_1\om_1\wedge \om_4+2vt_2 \om_2\wedge
\om_4-\a_{43}\wedge
\om_3 \label{eq30}\\
d(\a_{43})&=4v^2(t_1^2+t_2^2+1)\om_4\wedge\om_3\notag\\
dr&=-t_1(3r^2-rv+6v^2)\om_1+t_2(r-3v)(3r-2v)\om_2 \notag\\
dv&=-t_1v(7v+r)\om_1+t_2v(r-3v)\om_2 \notag\\
d(t_1)&=[r(t_1^2+t_2^2+1)+v(5t_1^2-3t_2^2+5)]\om_1\notag \\
d(t_2)&=8vt_1t_2\om_1+(v-r)(t_1^2+t_2^2+1)\om_2\notag
\end{align}
for some functions $t_1,t_2$. Differentiating these equations
yields only identities and the solution depends on 2 parameters
$t_1,t_2$.

We were not able to integrate completely the structure equations
and find the family of special Lagrangian 4-folds that are
solutions to these equations. The only thing we could observe is
that the generic solution has rank 2, since the following
relations hold between the parameters $r,v,t_1$ and $t_2$:
$$d\biggl(\frac{(t_1^2+t_2^2+1)v^{\frac{4}{5}}(r-3v)}{(r-v)^{\frac{3}{5}}}\biggl)=0 \ \mbox{  and  }
\
d\biggl(\frac{(t_2^2(r-3v)+(r-v)(t_2^2+1))v^{\frac{7}{5}}}{(r-v)^{\frac{4}{5}}}\biggl)=0$$
{\it Remark:} In principle, the structure equations can be
integrated using the reduction process for special Lagrangian
submanifolds with symmetries, by solving the ODE associated to it,
as in \cite{dj4}. The solution would be in terms of the Jacobi
elliptic functions. As for now, we do not have an explicit
integral yet.

\subsubsection{$SO(2)\ltimes S_3$-symmetry}

\begin{theorem}
Suppose that $L\subset \C^4$ is a connected special Lagrangian
$4$-fold with the property that its fundamental cubic at each
point has an $SO(2)\ltimes S_3$-symmetry. Then $L$ is congruent to
a product $\Sigma\ti \R^2$, where $\Sigma\subset \C^2$ is a
holomorphic curve.
\end{theorem}

\pf Let $L$ be a special Lagrangian 4-fold that satisfies the
hypotheses of the theorem and let $C$ be its fundamental cubic.
Proposition \ref{cont} implies that there exists a function
$r\colon L\to \R_+$ for which the equation
$$ C=r(\om_1^3-3\om_1\om_2^2)$$
defines an $SO(2)\ltimes S_3$-subbundle $F\subset P_L$ of the
$L$-adapted coframe bundle $P_L \ra L$, subbundle on which the
1-forms $\om_1,\om_2,\om_3,\om_4$ and $\alpha_{43}$ form a basis.
Similar calculations as in previous cases show that the structure
equations on $F$ are:
\begin{align}
d\om_1&=t_1\om_1\wedge\om_2, \quad d\om_2=t_2\om_1\wedge\om_2,
\quad
d\om_3=\a_{43}\wedge\om_4, \quad d\om_4=-\a_{43}\wedge\om_3 \notag\\
d\a_{43}&=0 \notag \\
dr&=-3rt_2\om_1+3rt_1\om_2 \label{eq20}\\
dt_1&=-u_2\om_1+(t_1^2+t_2^2-2r^2+u_1)\om_2 \notag \\
dt_2&=u_1\om_1+u_2\om_2\notag
\end{align}
for some functions $u_1,u_2$. Differentiation of these equations
does not lead to new relations among the quantities because the
system becomes involutive, according to Cartan-K{\"a}hler Theorem.
This is seen by computing the Cartan characters: $s_1=2$,
$s_2=s_3=s_4=0$ and noticing that the space of integral elements
at each point is parametrized by 2 parameters $u_1,u_2$.

 We are looking to integrate the above equations
and find the family of special Lagrangian 4-folds that satisfy the
hypothesis of the theorem. The Cartan-K{\"a}hler analysis tells us
that the solution should depend on 2 functions of one variable.

From the above structure equations, we can see that
$\om_1=\om_2=0$ and $\om_3=\om_4=0$ define integrable 2-plane
fields on $L$. The 2-dimensional leaves of the 2-plane field
$\Gamma_1$ defined by $\om_3=\om_4=0$ are congruent along
$\Gamma_2$, the codimension 2 foliation defined by
$\om_1=\om_2=0$. This is clear since $dt_1=dt_2\equiv 0 \mod
\{\om_1,\om_2\}$ and therefore the structure equations of
$\Gamma_1$ are:
$$d\om_1=t_1\om_1\wedge\om_2,\ d\om_2=t_2\om_1\wedge\om_2 $$
where $t_1, t_2$ are constant along $\Gamma_2$. Also, the third to
fifth equations in (\ref{eq20}) imply that the leaves of the
foliation $\Gamma_2$ are 2-planes which are congruent along
$\Gamma_1$ since $d(e_3\wedge e_4)=0$ and the 2-planes are real,
spanned by $\{e_3,e_4\}$.  Therefore, $L$ is a product
$L=\Sigma\ti \R^2$ where $\Sigma \subset \C^2$ is a surface. In
order for $L$ to be a special Lagrangian 4-fold, $\Sigma $ should
be a holomorphic curve with respect to a certain unique complex
structure on $\C^2$. This is because of the following argument.
Choose coordinates $z_k=x_k+iy_k,\ k=1...4$ on $L=\Sigma\ti \R^2$.
Then $L$ is special Lagrangian if and only if the 2-forms
$dx_1\wedge dy_1+dx_2\wedge dy_2$ and $dx_1\wedge dy_2+dy_1\wedge
dx_2$ each vanish when pulled back to $\Sigma$. But
$$ (dx_1\wedge dy_1+dx_2\wedge
dy_2)+i(dx_1\wedge dy_2+dy_1\wedge
dx_2)=(dx_1-idx_2)\wedge(dy_1+idy_2)=du\wedge dv$$ where
$u=x_1-ix_2$ and $v=y_1+iy_2$ are a different set of complex
coordinates on $\C^2$. Then $\Sigma \subset \C^2$ is special
Lagrangian if and only if $du\wedge dv\mid_{\Sigma}=0$, which says
that $\Sigma$ is a holomorphic curve in $\C^2$ with respect to the
complex coordinates $(u,v)$ on $\C^2$. $\Box$

\subsubsection{$SO(2)$-symmetry}
\begin{theorem}
Suppose that $L\subset \C^4$ is a connected special Lagrangian
$4$-fold with the property that its fundamental cubic has an
$SO(2)$-symmetry at each point. Then $L$ is invariant under an
$SO(3)$-action, whose orbits are $2$-spheres, and the surface we
obtain in the quotient by this action is a pseudo-holomorphic
curve, with respect to an almost complex structure.
\end{theorem}

\pf Let $L$ be a special Lagrangian 4-fold that satisfies the
hypotheses of the theorem and let $C$ be its fundamental cubic.
Proposition \ref{cont} implies that there exist functions $r\colon
L\to \R,\ s\colon L\to R_+$ and $v\colon L \to \R_+$ for which the
equation
$$ C=r(\om_1^3-3\om_1\om_2^2)+s(3\om_1^2\om_2-\om_2^3)+3v\om_1(\om_1^2+\om_2^2-2\om_3^2-2\om_4^2)$$
defines an $SO(2)$-subbundle $F\subset P_L$ of the $L$-adapted
coframe bundle $P_L \ra L$. The 1-forms $\om_1, \om_2, \om_3,
\om_4$ and $\alpha_{43}$ form a basis and satisfy the structure
equations:
\begin{align}
d\om_1&=[-t_1s+(r-3v)t_2]\om_1\wedge\om_2, \quad
d\om_2=[t_2s+(r-v)t_1]\om_1\wedge\om_2 \notag\\
d\om_3&=2t_1v\om_1\wedge\om_3+2t_2v\om_2\wedge\om_3+\a_{43}\wedge\om_4,
\quad
d\om_4=2t_1v\om_1\wedge\om_4+2t_2v\om_2\wedge\om_4-\a_{43}\wedge\om_3\notag\\
d\a_{43}&=-4v^2(t_1^2+t_2^2+1)\om_3\wedge\om_4 \notag\\
dr&=[t_1(-6v^2+vr-3s^2-3r^2)-11t_2vs-t_4]\om_1+\notag\\
+&[-t_1vs+t_2(6v^2-11vr+3r^2+3s^2)+t_3]\om_2
\notag \\
ds&=t_3\om_1+t_4\om_2 \label{eq73} \\
dv&=-v[t_1(r+7v)+t_2s]\om_1+v[-t_1s+(r-3v)t_2]\om_2 \notag\\
dt_1&=[t_1^2(r+5v)+t_2^2(r-3v)+r+5v]\om_1
+[s(t_1^2+t_2^2)+s]\om_2 \notag \\
dt_2&=[8vt_1t_2+s+s(t_1^2+t_2^2)]\om_1+(v-r)[(t_1^2+t_2^2)+1]\om_2 \notag\\
dt_3&=u_1\om_1+[u_2+t_2t_3(r-3v)+t_1t_4(r-v)-t_1t_3s+t_2t_4s]\om_2\notag
\\
dt_4&=u_2\om_1+(-u_1-9t_1^2s^3-6s^3+24rsv-6r^2s-18v^2s-60t_1t_2vs^2+90rsvt_2^2
\notag \\
+&3vt_1t_3+30t_1^2vrs-9t_2^2s^3-7rt_1t_3-9t_2^2r^2s-21t_1^2v^2s-7t_2t_3s+7t_2t_4r
\notag \\
-&25t_2t_4v-141t_2^2v^2s-9t_1^2r^2s-7t_1t_4s)\om_2 \notag
\end{align}
for some functions $u_1,u_2$. The above system is in involution,
so differentiation of these equations does not lead to new
relationships among quantities.

From the above structure equations, we can see that
$\om_1=\om_2=0$ and $\om_3=\om_4=0$ define integrable 2-plane
fields on $L$. The 2-dimensional leaves of the 2-plane field
$\Gamma_1$ defined by $\om_1=\om_2=0$ are 2-spheres. This is clear
since the structure equations of the leaves of $\Gamma_1$ are the
structure equations of a 2-dimensional sphere of constant radius
$4v^2(t_1^2+t_2^2+1)$:
$$d\om_3=-\a_{34}\wedge\om_4,\ d\om_4=\a_{34}\wedge\om_3, \ d\a_{34}=4v^2(t_1^2+t_2^2+1)\om_3\wedge\om_4 $$
and $t_1, t_2, v$ are constant along $\Gamma_1$ since $dt_1\equiv
dt_2\equiv dv \equiv 0 \mod \{\om_1,\om_2\}$. Therefore $L$ is
foliated by non-congruent spheres.

The 2-dimensional leaves of the other foliation $\Gamma_2$,
defined by $\om_3=\om_4=0$, are congruent. This follows from the
structure equations
\begin{align*}
d\om_1&=[-t_1s+(r-3v)t_2]\om_1\wedge \om_2, \quad
d\om_2=[t_2s+(r-v)t_1]\om_1\wedge \om_2 \notag
\end{align*}
and the fact that $dr\equiv dv\equiv dt_1\equiv dt_2 \equiv 0 \mod
\{\om_1,\om_2\}$.

Also, the structure equations imply $d(e_1\wedge e_2 \wedge Je_1
\wedge Je_2)=0 \mod \{\om_3,\om_4\}$. Therefore the complex
2-plane $(e_1,e_2,Je_1,Je_2)$ is constant along each leaf of the
$\Gamma_2$-foliation and each such leaf lies in an affine plane
parallel to this 2-plane.

If we let $\om_{21}=[-t_1s+(r-3v)t_2]\om_1+[t_2s+(r-v)t_1]\om_2$,
the structure equation for the $\Gamma_2$ leaves can be written
as:
$$d\om_1=\om_{21}\wedge \om_2, \ d\om_2=-\om_{21}\wedge \om_1, \
d\om_{21}=2(r^2+s^2-v^2)\om_1 \wedge \om_2$$ This shows that the
leaves of the $\Gamma_2$ foliation are congruent surfaces of Gauss
curvature $2(v^2-r^2-s^2)$, lying in the affine complex 2-plane
$(e_1,e_2,Je_1,Je_2)$.

Computations show that the structure equations are invariant under
an $SO(3)$-rotation about some point in $\C^4$. Therefore, the
solutions should be special Lagrangian 4-folds that are invariant
under the subgroup $SO(3)$, as it sits naturally in $\so$ and
hence in $SU(4)$. The orbits of the $SO(3)$-action are 2-spheres.

We look now for special Lagrangian 4-folds $L$, invariant under
the action of $SO(3)$. Let $${\bf z}=\begin{pmatrix}{\bf x}+i{\bf
y}\\x_4+iy_4\end{pmatrix}, \ {\bf x}=(x_1,x_2,x_3), {\bf
y}=(y_1,y_2,y_3) \in\R^3$$ denote the coordinates on $\C^4$. The
subgroup $SO(3)$ acts diagonally by rotation in ${\bf x}$ and
${\bf y}$,
$$A\cdot \begin{pmatrix}{\bf x}+i{\bf y}\\x_4+iy_4\end{pmatrix}=
\begin{pmatrix}A{\bf x}+iA{\bf y}\\x_4+iy_4\end{pmatrix},\ A\in SO(3)$$

Let $X_1,X_2,X_3$ be the infinitesimal generators of $SO(3)$,
where
\begin{align*}
X_1=&x_2\frac{\partial}{\partial x_3}-x_3\frac{\partial}{\partial
x_2}+y_2\frac{\partial}{\partial y_3}-y_3\frac{\partial}{\partial y_2}\\
X_2=&x_3\frac{\partial}{\partial x_1}-x_1\frac{\partial}{\partial
x_3}+y_3\frac{\partial}{\partial y_1}-y_1\frac{\partial}{\partial y_3}\\
X_3=&x_1\frac{\partial}{\partial x_2}-x_2\frac{\partial}{\partial
x_1}+y_1\frac{\partial}{\partial y_2}-y_2\frac{\partial}{\partial
y_1}
\end{align*}
The 4-fold $L$ is invariant under the flow of $X_i$, $i=1,2,3$, so
$(X_i \lrcorner \om)\mid_L=0, i=1,2,3,$ where $\om=d{\bf x}\cdot
d{\bf y}+dx_4\wedge dy_4$ is the symplectic form and
$$d{\bf x}\cdot d{\bf y}:=dx_1\wedge
 dy_1+dx_2\wedge dy_2+dx_3\wedge dy_3$$

It is easy to calculate that $$(X_1 \lrcorner
\om)\mid_L=d(x_2y_3-x_3y_2)$$ which implies that
$x_2y_3-x_3y_2=c_1$, where $c_1\in\R$ is a constant. Similarly, we
can show that
$$x_3y_1-x_1y_3=c_2, \ x_1y_2-x_2y_1=c_3$$
From here it follows that $${\bf x}\times {\bf y}={\bf
c}=(c_1,c_2,c_3),$$ where ${\bf c}\in\R^3$ is a constant vector.

If $c\not =0$, then ${\bf x},{\bf y}$ are linearly independent and
therefore the stabilizer of a point on the orbit is trivial. This
implies that the orbit has dimension 3, but we know that the
orbits are 2 dimensional spheres. It follows that $c=0$, i.e.
${\bf x}$ and ${\bf y}$ are linearly dependent. So, $L$ lies in
the 6-manifold $\Sigma\subset\C^4$ on which the coordinates are
given by
$${\bf z}=\begin{pmatrix}(x+iy){\bf u}\\x_4+iy_4\end{pmatrix}, \ x,y\in\R, \ {\bf u}=(u_1,u_2,u_3)\in S^2$$

It is easy to compute that
\begin{align*}
\om\mid _\Sigma=&d(xu_1)\wedge d(yu_1)+d(xu_2)\wedge
d(yu_2)+d(xu_3)\wedge
d(yu_3)+dx_4\wedge dy_4\\
=&dx\wedge dy+dx_4\wedge dy_4
\end{align*}
Dividing out by the sphere action, we obtain in the quotient a
4-dimensional manifold $X^4$, with coordinates $(x,y,x_4,y_4)$.
The leaves of the $\om_3=\om_4=0$ foliation are 2-dimensional
surfaces $M^2$.
 We calculate now the pullback of the volume form
to $\Sigma$. Denote $z=x+iy$ and $w=x_4+iy_4$. So,
\begin{align*}
\Omega\mid _\Sigma&=dz_1\wedge dz_2\wedge dz_3\wedge dz_4\mid_\Sigma =d(zu_1)\wedge d(zu_2)\wedge d(zu_3)\wedge dw\\
&=z^2(u_3du_1\wedge du_2+u_1du_2\wedge du_3+u_2du_3\wedge
du_1)\wedge dz \wedge dw=\frac{1}{3}d(z^3)\wedge dw \wedge d\sigma
\end{align*}
 where $d\sigma=u_3du_1\wedge du_2+u_1du_2\wedge du_3+u_2du_3\wedge du_1$ is the area
 form of the 2-sphere $S^2$.

Then $L\subset \Sigma$ is special Lagrangian if and only if the
2-forms
\begin{align*}
\a&=dx\wedge dy+dx_4\wedge dy_4=\frac{i}{2}(dz\wedge d{\bar z}+dw\wedge d{\bar w})\\
\b&=\Im (\frac{1}{3}d(z^3)\wedge dw)=-\frac{i}{6}(d(z^3)\wedge
dw-d({\bar z^3})\wedge d{\bar w})
\end{align*}
each vanish when pulled back to $M^2\subset X^4$. But:
\begin{align*}
\a \wedge \a&=-\frac{1}{2}(dz\wedge d{\bar z}\wedge dw\wedge
d{\bar w})\\
\b\wedge \b&=\frac{1}{18}(d(z^3)\wedge dw\wedge d({\bar
z^3})\wedge d{\bar w})=\frac{1}{2}((z{\bar z})^2dz\wedge dw\wedge
d{\bar z}\wedge d{\bar w})
\end{align*}
Rescaling $\b$ by dividing it by $z{\bar z}$
$$\tilde{\b}=-\frac{i}{6z{\bar z}}(d(z^3)\wedge dw-d({\bar z^3})\wedge d{\bar w})$$
we get that $(\a+i\tilde{\b})^2=0$, so this form is decomposable
and it is easy to compute that
$$ \a+i\tilde{\b}=\frac{i}{2z{\bar z}}(\xi_1\wedge \xi_2),$$
where $\xi_1=zdz-i{\bar z}d{\bar w}$ and $\xi_2=zd{\bar z}-izdw $.
The forms $\xi_1,\xi_2$ form a system which is not integrable
since
$$d\xi_1=\frac{z}{{\bar z}}dw\wedge d{\bar w}\not =0 \mod \{\xi_1,\xi_2\}$$
In fact there is no combination of the forms $\xi_1$, $\xi_2$ that
is integrable.

 Since $\a+i\tilde{\b}\mid_M=0$, it implies that $M^2$ is a
pseudo-holomorphic curve in $X^4$, with respect to a certain
almost complex structure, which is not integrable. Conversely,
every almost complex surface in $X^4$ lifts to a special
Lagrangian 4-fold $L\subset \Sigma^6\subset \C^4$. $\Box$

%%%%%%%%%%%%%%%%%%%%%%%%%%%%%%%%%%%%%%%%%%
%%%%%%Discrete symmetry
%%%%%%%%%%%%%%%%%%%%%%%%%%%%%%%%%%%%%%%%%
 \subsection{Discrete symmetry}

Next, we are analyzing the case of discrete symmetry. Suppose that
the stabilizer $G$ of the fundamental cubic of a special
 Lagrangian 4-fold is a finite subgroup of $\so$. If $g$ is an element of $G$,
 then $g$ is conjugate to an
 element in the maximal torus of $\so$:
$$\biggl\{\begin{pmatrix}e^{2\pi ir} &0\\ 0& e^{2\pi
is}\end{pmatrix}, r\in\Q, \ s\in \Q, \ r,s<1 \biggl\}$$

The following result tells us when there exists a harmonic cubic
in 4 variables fixed by a nontrivial element $g$ in the maximal
torus.

\begin{proposition}
\label{discrete}
 The element $g=\begin{pmatrix}e^{2\pi ir} &0\\ 0& e^{2\pi
is}\end{pmatrix}$, where $r\in\Q, s\in \Q, \ r,s<1$, fixes a
nontrivial harmonic cubic in four variables $(x_1,x_2,x_3,x_4)$ if
and only if at least one of the following conditions is satisfied:
\\

$\quad 1.\ 3r\in\Z, \mbox{when the fixed harmonic cubics contain
linear combinations of }$ \vskip.2cm
 $\quad \{x_1^3-3x_1x_2^2,\
3x_1^2x_2-x_2^3\};$ \vskip.2cm

 $\quad 2.\ r\in \Z, \mbox{when the
fixed harmonic cubics contain linear combinations of }$ \vskip.2cm
 $\quad \{x_1^3-3x_1x_2^2,\ 3x_1^2x_2-x_2^3,\
  x_1(x_1^2+x_2^2-2x_3^2-2x_4^2),\
x_2(x_1^2+x_2^2-2x_3^2-2x_4^2)\}; $\vskip.2cm

 $\quad 3.\
2r+s\in\Z, \mbox{when the fixed harmonic cubics contain linear
combinations of}$ \vskip.2cm
 $\quad \{(x_1^2-x_2^2)x_3-2x_1x_2x_4,
\ (x_1^2-x_2^2)x_4+2x_1x_2x_3\};$\vskip.2cm

 $\quad 4.\ 2r-s\in\Z,
\mbox{when the fixed harmonic cubics contain linear combinations
of}$\vskip.2cm
 $\quad \{(x_1^2-x_2^2)x_3+2x_1x_2x_4, \
  (x_1^2-x_2^2)x_4-2x_1x_2x_3\};$\vskip.2cm

$\quad 5.\ 2s+r\in\Z, \mbox{when the fixed harmonic cubics contain
linear combinations of }$\vskip.2cm
 $\quad \{(x_3^2-x_4^2)x_1-2x_2x_3x_4,
\ (x_3^2-x_4^2)x_2+2x_1x_3x_4\};$\vskip.2cm

 $\quad 6.\ 2s-r\in\Z,
\mbox{when the fixed harmonic cubics contain linear combinations
of}$\vskip.2cm $\quad \{(x_3^2-x_4^2)x_1+2x_2x_3x_4, \
  (x_3^2-x_4^2)x_2-2x_1x_3x_4\};$\vskip.2cm

$\quad 7.\ 3s\in\Z, \mbox{when the fixed harmonic cubics contain
linear combinations of}$ \vskip.2cm
 $\quad  \{x_3^3-3x_3x_4^2,\
3x_3^2x_4-x_4^3\};$\vskip.2cm

 $\quad 8.\ s\in\Z, \mbox{when the fixed
harmonic cubics contain linear combinations of}$\vskip.2cm
 $\quad \{x_3^3-3x_3x_4^2,\ 3x_3^2x_4-x_4^3,\
  x_3(2x_1^2+2x_2^2-x_3^2-x_4^2),\ x_4(2x_1^2+2x_2^2-x_3^2-x_4^2)\}$.
\end{proposition}

\pf Let ${\mathcal P}^3_\C={\mathcal
P}^3(z_1,z_2,\bar{z_1},\bar{z_2})$ be the space of complexified
cubic polynomials, in the variables $(z_1,z_2,{\bar z_1},{\bar
z_2})$, where $z_1=x_1+ix_2,\ z_2=x_3+ix_4$. The maximal torus in
$SO(4)$ acts on ${\mathcal H}^3_\C={\mathcal
H}^3(z_1,z_2,\bar{z_1},\bar{z_2})$, the space of complexified
harmonic cubics in 4 variables $(z_1,z_2,{\bar z_1},{\bar z_2})$,
as follows:
$$\bigl(\begin{smallmatrix}e^{2\pi ir} &0\\ 0& e^{2\pi
is}\end{smallmatrix}\bigl).P(z_1,z_2,\bar{z_1},\bar{z_2})=P(z_1^*,z_2^*,\bar{z_1^*},\bar{z_2^*}),
\ P\in {\mathcal H}^3_\C$$ where $z_1^*=e^{2\pi ir}z_1,\
z_2^*=e^{2\pi is}z_2$. Under this action, the space ${\mathcal
H}^3_\C$ decomposes as follows:
$${\mathcal H}^3_\C={\mathcal H}^3(z_1,z_2)\oplus
{\mathcal H}({\mathcal P}^2(z_1,z_2)\otimes {\mathcal
P}^1(\bar{z_1},\bar{z_2}))\oplus{\mathcal H}({\mathcal
P}^1(z_1,z_2)\otimes {\mathcal
P}^2(\bar{z_1},\bar{z_2}))\oplus{\mathcal
H}^3(\bar{z_1},\bar{z_2})$$

A basis for the space ${\mathcal H}^3(z_1,z_2)$ is given by the
polynomials $\{z_1^3, z_1^2z_2, z_1z_2^2, z_2^3\}$. Since
$g.z_1=e^{2\pi ir} z_1$ and $g.z_2=e^{2\pi is} z_2$, it follows
that $g.z_1^3=e^{6\pi ir}z_1^3$, $g.z_1^2z_2=e^{2\pi
i(2r+s)}z_1^2z_2$, $g.z_1z_2^2=e^{2\pi(r+2s)}z_1z_2^2$ and
$g.z_2^3=e^{6\pi is}z_2^3$. This further implies that unless
$e^{6\pi ir}=1$, $e^{2\pi i(2r+s)}=1$, $e^{2\pi i(r+2s)}=1$ or
$e^{6\pi is}=1$, there is no fixed element in the space ${\mathcal
H}^3(z_1,z_2)$. The above conditions are equivalent to $3r\in\Z$,
$2r+s\in\Z$, $r+2s\in\Z$ or $3s\in\Z$. Therefore, ${\mathcal
H}^3(z_1,z_2)$ decomposes into the following four weight spaces:
$${\mathcal H}^3(z_1,z_2)=V_{(3,0)}\oplus V_{(2,1)}\oplus V_{(1,2)} \oplus V_{(0,3)}$$
All these weight spaces have multiplicity 1 and a basis in
$V_{(3,0)}, V_{(2,1)}, V_{(1,2)}, V_{(0,3)}$ is given by the
harmonic polynomials $z_1^3,z_1^2z_2,z_1z_2^2$ and $z_2^3$
respectively.

We analyze now the fixed elements for the space ${\mathcal
H}({\mathcal P}^2(z_1,z_2)\otimes {\mathcal
P}^1(\bar{z_1},\bar{z_2}))$ of harmonic polynomials in ${\mathcal
P}^2(z_1,z_2)\otimes {\mathcal P}^1(\bar{z_1},\bar{z_2})$. It is
easy to see that a basis in the space ${\mathcal H}({\mathcal
P}^2(z_1,z_2)\otimes {\mathcal P}^1(\bar{z_1},\bar{z_2}))$ is
given by the harmonic cubics $$\{z_1^2{\bar z_2},\ z_1^2{\bar
z_1}-2z_1z_2{\bar z_2},\ z_2^2{\bar z_2}-2z_1z_2{\bar z_1},\
z_2^2{\bar z_1}\}$$ and this space decomposes into:
$${\mathcal H}({\mathcal
P}^2(z_1,z_2)\otimes{\mathcal
P}^1(\bar{z_1},\bar{z_2}))=V_{(2,-1)}\oplus V_{(1,0)}\oplus
V_{(0,1)} \oplus V_{(-1,2)}$$ We can see that unless at least one
of the conditions: $2r-s\in \Z, r\in\Z, s\in\Z$ or $2s-r\in \Z$
are satisfied, there is no fixed vector in any of the weight
spaces.

 Doing a similar argument, one can see that a basis in the space ${\mathcal H}({\mathcal P}^1(z_1,z_2)\otimes {\mathcal
P}^2(\bar{z_1},\bar{z_2}))$ is given by the polynomials
$$\{z_2{\bar z_1}^2,\ z_1{\bar z_1}^2-2z_2{\bar z_1}{\bar z_2},\
z_2{\bar z_2}^2-2z_1{\bar z_1}{\bar z_2},\ z_1{\bar z_2}^2\}$$ and
the space decomposes into:
$${\mathcal H}({\mathcal
P}^1(z_1,z_2)\otimes{\mathcal
P}^2(\bar{z_1},\bar{z_2}))=V_{(-2,1)}\oplus V_{(-1,0)}\oplus
V_{(0,-1)} \oplus V_{(1,-2)}$$ Unless at least one of the
conditions: $-2r+s\in \Z, -r\in\Z, -s\in\Z$ or $-2s+r\in \Z$ is
satisfied, there is no fixed vector in any of the weight spaces.

Finally, a basis in the space ${\mathcal H}^3({\bar z_1},{\bar
z_2})$ is given by the polynomials $\{{\bar z_1}^3, {\bar
z_1}^2{\bar z_2}, {\bar z_1}{\bar z_2^2}, {\bar z_2}^3\}$ and this
space decomposes into the weight spaces:
$${\mathcal H}^3({\bar z_1},{\bar z_2})=V_{(-3,0)}\oplus V_{(-2,-1)}\oplus V_{(-1,-2)} \oplus V_{(-0,-3)}$$
For there to be a fixed vector in this space, at least one of the
following conditions should be satisfied: $-3r\in\Z$,
$-2r-s\in\Z$, $-r-2s\in\Z$ or $-3s\in\Z$.

Therefore, the space of complexified harmonic cubics decomposes
under the action of the maximal torus into 8 pairs of opposite
weight spaces, each of multiplicity one:
\begin{align*}
{\mathcal H}^3_\C&=V_{(3,0)}\oplus V_{(-3,0)}\oplus
V_{(2,1)}\oplus V_{(-2,-1)}\oplus
V_{(1,2)}\oplus V_{(-1,-2)}\oplus V_{(0,3)}\oplus V_{(0,-3)}\oplus V_{(2,-1)}\\
&\oplus V_{(-2,1)}\oplus V_{(1,0)}\oplus V_{(-1,0)}\oplus
V_{(0,1)}\oplus V_{(0,-1)}\oplus V_{(1,-2)}\oplus V_{(-1,2)}
\end{align*}

A real harmonic cubic is the sum of elements drawn from these
weight spaces, with the coefficients in opposite weight spaces
being complex conjugates. Then, there exists a fixed element in
the space of real harmonic cubics in 4 variables if and only if
there are nontrivial elements in the maximal torus that act
trivially on at least one pair of these weight spaces. By the
above analysis, one can see that this is equivalent to the
satisfaction of at least one of the following conditions: (1)
$3r\in\Z$, (2) $r\in\Z$, (3) $2r+s\in\Z$, (4) $2r-s\in\Z$, (5)
$2s+r\in\Z$, (6) $2s-r\in\Z$, (7) $3s\in\Z$, (8) $s\in\Z$. Next we
assume that exactly one of the conditions above is satisfied:

1) $3r\in \Z$. In this case $g$ acts trivially on the pair
 of opposite weight spaces $V_{(3,0)}$ and $V_{(-3,0)}$ and the fixed real harmonic cubics in 4 variables
  are of the form $az_1^3+{\bar a}{\bar z_1}^3$. So,
 $$C={\mbox Re}(az_1^3),$$ where $a\in\C$. Therefore,
 a basis in the space of fixed real harmonic cubics is given by the
 harmonic polynomials $\{x_1^3-3x_1x_2^2,\ 3x_1^2x_2-x_2^3\}$.

2) $r\in \Z$. This condition implies also $3r\in \Z$ and $g$ acts
trivially on the pairs of
  opposite weight spaces $V_{(3,0)}, V_{(-3,0)}, V_{(1,0)}$ and
  $V_{(-1,0)}$. So, the fixed real harmonic cubics are of the form:
  $$C={\mbox Re}(az_1^3+b(z_1^2{\bar z_1}-2z_1z_2{\bar z_2})),$$ where $a, b\in \C$. Therefore,
 a basis in the space of fixed real harmonic cubics is: $$\{x_1^3-3x_1x_2^2,\ 3x_1^2x_2-x_2^3, \ x_1(x_1^2+x_2^2-2x_3^2-2x_4^2),\
 x_2(x_1^2+x_2^2-2x_3^2-2x_4^2)\}$$

3) $2r+s \in \Z$. Then $g$ acts trivially on the pair
 of opposite weight spaces $V_{(2,1)}$ and $V_{(-2,-1)}$ and the fixed real
 harmonic cubics are of the form:
 $$C={\mbox Re}(az_1^2z_2),$$ where $a\in\C$. A basis in the space of fixed real harmonic cubics is given by the
polynomials: $\{(x_1^2-x_2^2)x_3-2x_1x_2x_4, \
 (x_1^2-x_2^2)x_4+2x_1x_2x_3\}$

4) $2r-s\in \Z$. The element $g$ acts trivially on $V_{(2,-1)}$
and $V_{(-2,1)}$ and the fixed real
 harmonic cubics are of the form:
 $$C={\mbox Re}(az_1^2{\bar z_2}),$$ where $a\in\C$. Thus,
 a basis in the space of fixed real harmonic cubics is given by the polynomials: $\{(x_1^2-x_2^2)x_3+2x_1x_2x_4, \
  (x_1^2-x_2^2)x_4-2x_1x_2x_3\}$

5) $2s+r\in \Z$. In this case $g$ acts trivially $V_{(1,2)}$ and
$V_{(-1,-2)}$ and $$C={\mbox Re}(az_1z_2^2),$$ where $a\in\C$ is
the general harmonic cubic polynomial fixed by the action.
Therefore,
 a basis for the space of fixed real harmonic cubics is given by the
 harmonic polynomials $\{(x_3^2-x_4^2)x_1-2x_2x_3x_4, \
  (x_3^2-x_4^2)x_2+2x_1x_3x_4\}$.

6) $2s-r\in \Z$. Then $g$ acts trivially on the pair
 of opposite weight spaces $V_{(-1,2)}$ and $V_{(1,-2)}$ and the fixed real
 harmonic cubics are of the form:
 $$C={\mbox Re}(z_2^2{\bar z_1}),$$ where $a\in\C$. A basis is given by the polynomials: $\{(x_3^2-x_4^2)x_1+2x_2x_3x_4, \
  (x_3^2-x_4^2)x_2-2x_1x_3x_4\}$

7) $3s\in \Z$. Now $g$ acts trivially on the pair
 of opposite weight spaces $V_{(0,3)}$ and $V_{(0,-3)}$ and the fixed real
 harmonic cubics are of the form:
 $$C={\mbox Re}(az_2^3),$$ where $a\in\C$ and therefore, a basis is given by $\{x_3^3-3x_3x_4^2,\ 3x_3^2x_4-x_4^3\}$.

8) $s\in \Z$. In this last case, $g$ acts trivially on $V_{(0,3)},
V_{(0,-3)}, V_{(0,1)}$ and
  $V_{(0,-1)}$. The real harmonic cubics fixed by the action are of the form:
  $$C={\mbox Re}(az_2^3+b(z_2^2{\bar z_2}-2z_1z_2{\bar z_1})),$$ where
$a, b\in \C$. Therefore, a basis in the space of fixed real
harmonic cubics is: $$\{x_3^3-3x_3x_4^2,\ 3x_3^2x_4-x_4^3, \
x_3(2x_1^2+2x_2^2-x_3^1-x_4^2),\
 x_4(2x_1^2+2x_2^2-x_3^2-x_4^2)\}  \quad \Box$$
\vskip.5cm

{\it Remark 1:} In Figure 1 below we graphed in the coordinates
$(r,s)$ mod $\Z$ all the possibilities appearing in Proposition
\ref{discrete}. By moding out by the Weyl group of $\so$, we can
consider the possibilities only in the triangle found by
intersecting the regions below the lines $r=s$ and $s=1-r$.
Furthermore, since the stabilizer of $g$ in $\so$ coincides with
its stabilizer in $O(4)$, we can actually mod out by the Weyl
group of $O(4)$. The elements $\l(
\begin{smallmatrix} e^{2\pi ir}&0\\0&e^{2\pi is}\end{smallmatrix} \r)$ and $\l(
\begin{smallmatrix} e^{2\pi ir}&0\\0&e^{-{2\pi is}}\end{smallmatrix}\r)$ are conjugate
to each other in $O(4)$, by the element $\l( \begin{smallmatrix}
  1&0&0&0\\0&1&0&0\\0&0&-1&0\\0&0&0&1\end{smallmatrix}\r)\in O(4)$.
Therefore, we can restrict our attention to the cases found in the
small shaded triangle shown in Figure 1, which is a fundamental
Weyl chamber. \\

%\begin{figure}[htbp]
%\begin{center}
%\includegraphics[width=8.5cm,height=8.5cm]{figura.eps}
%\caption{The weight spaces and a fundamental Weyl chamber}
%\label{figura}
%\end{center}
%\end{figure}
\centerline{\psfig{figure=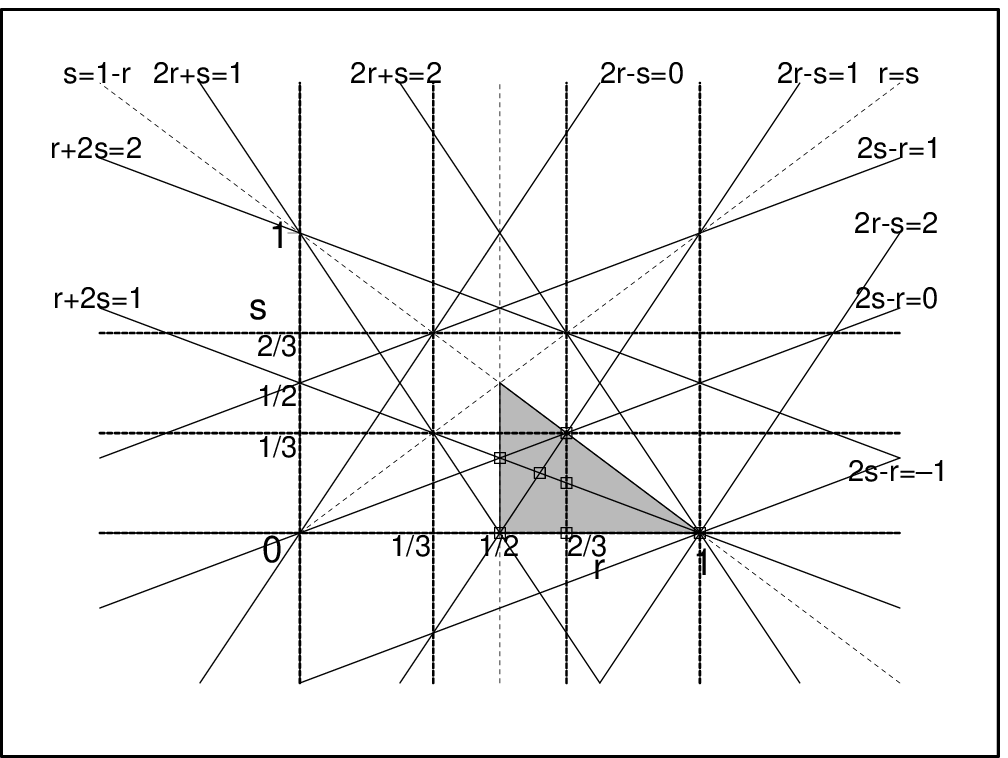,width=9.5cm, height=9.5cm}}
\begin{center}{Figure 1: The weight spaces and a fundamental Weyl chamber}\end{center}

\non{\it Remark 2:} If only one of the conditions in the small
triangle are satisfied, the stabilizer $G$ will be at least an
$SO(2)$, thus continuous, and we will recover the cases already
studied in section 3.2. For example, if only $2s-r\in\Z$, moding
out by $\Z$, we get $2s-r=0$. The stabilizer group in this case
looks like $\l\{\l(
\begin{smallmatrix}e^{2i\theta}&0\\0&e^{i\theta}\end{smallmatrix}\r),\theta\in\R\r\}$
and this is the case of continuous symmetry studied in Theorem
3.2.3.

Therefore, in order to find the fundamental cubics that have
discrete nontrivial stabilizers under the action of $\so$, we have
to look at elements that have nontrivial components in at least
two non-opposite weight spaces. As seen in Figure 1, up to
conjugacy in $O(4)$, there are six nontrivial elements in the
maximal torus that act trivially on more than two pairs of weight
spaces.
\begin{cor}
\label{cor1}
 If $G$ is a nontrivial discrete subgroup of $\so$ that stabilizes a nontrivial
polynomial $h\in {\mathcal H}_3(\R^4)$, then $G$ can not have
elements of order $>6$.
\end{cor}
\pf This follows from Proposition \ref{discrete} and the above
remarks. Any
element in $G$ is conjugate to an element of the form $g=\l(\begin{smallmatrix}e^{2\pi ir} &0\\
0& e^{2\pi is}\end{smallmatrix} \r)$, where $r=\frac{m}{p}\in\Q$,
\ $s=\frac{n}{q}\in \Q, \; r,s<1$.

By looking at the small triangle in Figure 1, we can see that
there are the following possibilities for the values of $r$ and $s
\mod \Z$ and mod the Weyl group:

\vskip.2cm

1) If $r+2s\in\Z$ and $3r\in\Z$, then $r=\frac{2}{3}$ and
$s=\frac{1}{6}$. The element $g=\l( \begin{smallmatrix}
  e^{\frac{4\pi i}{3}}&0\\ 0& e^{\frac{\pi i}{3}}
  \end{smallmatrix}\r)$ has {\bf order 6} and acts trivially on the
  pairs of opposite
  weight spaces $V_{(3,0)}, V_{(-3,0)}$ and $V_{(1,2)}, V_{(-1,-2)}$. The general harmonic cubic
  stabilized by this element is
  $$C={\mbox
  Re}(az_1^3+bz_1z_2^2), \ a,b\in\C$$
2) If $r+2s\in\Z$ and $2r-s\in\Z$, then we get
  $r=\frac{3}{5}$ and $s=\frac{1}{5}$. The element $g=\Bigl( \begin{smallmatrix}
  e^{\frac{6\pi i}{5}}&0\\ 0& e^{\frac{2\pi i}{5}}
  \end{smallmatrix}\Bigl)$ has {\bf order 5} and acts trivially on the
  pairs of opposite
  weight spaces $V_{(1,2)}, V_{(-1,-2)}$ and $V_{(2,-1)}, V_{(-2,1)}$. Therefore, the general harmonic cubic
  stabilized by this element is
  $$C=\Re(az_1z_2^2+bz_1^2{\bar z_2}), \ a,b\in\C$$
 3) If $2s-r\in\Z$ and $r+2s\in\Z$, then we get
$(r,s)=(\frac{1}{2},\frac{1}{4})$. The element $g=\Bigl(
\begin{smallmatrix}
  -1&0\\ 0& i \end{smallmatrix}\Bigl)$ has
  {\bf order 4} and acts trivially on the pairs of opposite
  weight spaces $V_{(-1,2)}, V_{(1,-2)}$ and $V_{(1,2)}, V_{(-1,-2)}$. The general harmonic cubic
  stabilized by this element is
  $$C={\mbox Re}(a{\bar z_1}z_2^2+bz_1z_2^2), \ a,b\in\C$$
4) If $s\in\Z$ and $3r\in\Z$, then we get
  $(r,s)=(\frac{2}{3},0)$. The element $g=\bigl( \begin{smallmatrix}
  e^{\frac{4\pi i}{3}}&0\\ 0& I_2\end{smallmatrix}\bigl)$ has {\bf
  order 3} and acts trivially on the pairs of opposite
  weight spaces $V_{(3,0)}, V_{(-3,0)},V_{(0,1)},V_{(0,-1)}$ and
  $V_{(0,3)}, V_{(0,-3)}$. The general harmonic cubic fixed by this element is
  $$C=\Re(az_1^3+bz_2^3+c(z_2^2{\bar z_2}-2z_1{\bar z_1}z_2)), \
  a,b,c\in\C$$
5) If $2s-r\in\Z$, $2r-s\in\Z$, $3r\in\Z$ and $3s\in\Z$ then we
get
  $(r,s)=(\frac{2}{3},\frac{1}{3})$. The element $g=\Bigl( \begin{smallmatrix}
  e^{\frac{4\pi i}{3}}&0\\ 0& e^{\frac{2\pi i}{3}}
  \end{smallmatrix}\Bigl)$ has {\bf order 3} and acts trivially on the
  pairs of opposite
  weight spaces $V_{(-1,2)}, V_{(1,-2)},V_{(2,-1)},V_{(-2,1)},
  V_{(3,0)}, V_{(-3,0)}$ and $V_{(0,3)}, V_{(0,-3)}$. The general harmonic cubic
  stabilized by this element is
  $$C=\Re(a{\bar z_1}z_2^2+bz_1^2{\bar
  z_2}+cz_1^3+ez_2^3), \ a,b,c,e\in\C$$
6) If $s\in\Z,\ 2r+s\in \Z$ and $2r-s\in\Z$, then we get
  $(r,s)=(\frac{1}{2},0)$. The element $g=\bigl( \begin{smallmatrix}
  -I_2 &0\\ 0& I_2\end{smallmatrix}\bigl)$ has {\bf order 2} and acts
  trivially on the pairs of opposite
  weight spaces $V_{(0,1)}, V_{(0,-1)},V_{(2,1)},V_{(-2,-1)},\\
  V_{(2,-1)}, V_{(-2,1)}$ and $V_{(0,3)}, V_{(0,-3)}$. The general harmonic cubic
  stabilized by this element is
  $$C=\Re(az_2^3+b(z_2^2{\bar
  z_2}-2z_1{\bar z_1}z_2)+cz_1^2z_2+ez_1^2{\bar z_2}),\ a,b,c,e\in\C$$
7) If $r\in\Z$, $s\in\Z$, $r+2s\in\Z$, $2r+s\in\Z$, $2s-r\in\Z$
and
  $2r-s\in\Z$, meaning all the conditions are satisfied at once, then
  we get
  $r=1$ and $s=0$, so $g$ is just the {\bf identity} element.

\subsubsection{Polyhedral symmetry}

Now we are going to find the nontrivial harmonic cubic polynomials
in 4 variables whose stabilizer is one of the polyhedral subgroups
of $\so$ described in Section 3.1. and we will study the families
of special Lagrangian 4-folds with these symmetries.

\begin{proposition}
\label{polyhedral}
 The $\so$-stabilizer of $C\in {\mathcal H}_3(\R^4)$ is a polyhedral
subgroup of $\so$ if and only if $C$ lies on the $\so$-orbit of
exactly one of the following polynomials:
\begin{align*}
1.& \ rx_1(x_1^2-x_2^2-x_3^2-x_4^2)+sx_2x_3x_4, \mbox { for some
$r,s>0$
satisfying $s\not = 2\sqrt 5 r$, whose stabilizer}\\
& \mbox{ is the tetrahedral
subgroup $\T$ of $\so;$}\\
2.& \ sx_2x_3x_4, \mbox { for some $s>0$, whose stabilizer is the
irreducibly
   acting octahedral subgroup $\O^+$;}\\
3.& \ r[x_1(x_1^2-x_2^2-x_3^2-x_4^2)+2\sqrt5 x_2x_3x_4],\
   r>0, \mbox {whose stabilizer is the irreducibly acting}\\
&\mbox{ icosahedral subgroup $\I^+$.}
\end{align*}
\end{proposition}
\pf The polyhedral subgroups of $\so$ were found to be the
tetrahedral subgroup $\T$ of order 12, the reducibly and
irreducibly acting octahedral subgroups $\O$ and $\O^+$, each of
order 24, the reducibly and irreducibly acting icosahedral
subgroups $\I$ and $\I^+$, each of order 60.

First we look at the tetrahedral subgroup $\T$ and find the
harmonic cubics in 4 variables $\{x_1,x_2,x_3,x_4\}$ that are
stabilized by this subgroup. As we have seen in Section 3.1,
$\T=\{[t,t]|t\in{\bf T}\}$, where ${\bf T}=\{\pm {\bf 1}, \pm {\bf
i}, \pm {\bf j},  \pm {\bf k}, {1\over2}(\pm {\bf 1}\pm{\bf i}\pm
{\bf
 j}\pm{\bf k})\}$ is the binary tetrahedral subgroup of the unit
quaternion group $U$, of order 24. The subgroup $\T$ sits in
$SO(3)$ and it is generated by the transformations: $[{\bf i},{\bf
i}]$ with representing matrix
$\biggl[\begin{smallmatrix}1&0&0&0\\0&1&0&0\\0&0&-1&0\\0&0&0&-1
\end{smallmatrix}\biggl]$, relative to the basis $\{{\bf 1}, {\bf i}, {\bf j}, {\bf
k}\}$, $[{\bf j},{\bf j}]$ with representing matrix
$\biggl[\begin{smallmatrix}1 &0&0&0\\0&-1&0&0\\0&0&1&0\\0&0&0&-1
\end{smallmatrix}\biggl]$ and $[{1\over2}({\bf 1}+{\bf i}+ {\bf
 j}+{\bf k}),{1\over2}({\bf 1}+{\bf i}+ {\bf
 j}+{\bf k})]$ with
representing matrix
$\biggl[\begin{smallmatrix}1&0&0&0\\0&0&1&0\\0&0&0&1\\0&1&0&0
 \end{smallmatrix}\biggl]$. We can see that $\T$ fixes a cubic in
 $\{x_1,x_2,x_3,x_4\}$ if and only if the cubic is invariant under the flips of the signs of two of
 the coordinates $\{x_2,x_3,x_4\}$ and also under permuting
 $\{x_2,x_3,x_4\}$ while keeping $x_1$ fixed. Therefore, the cubic
 should be a linear combination of the polynomials
 $x_1^3, \ x_1(x_2^2+x_3^2+x_4^2)$ and $x_2x_3x_4$. Now, considering
 the extra condition that the cubic should be harmonic, it
 follows that the harmonic cubics stabilized by $\T$ lie on the
 $\so$-orbit of the polynomial
\begin{equation}
\label{polyhedrap}
 C=rx_1(x_1^2-x_2^2-x_3^2-x_4^2)+sx_2x_3x_4,
\end{equation}
 for some $r,s\geq 0$.

We turn now to the harmonic cubics invariant under the reducibly
acting octahedral subgroup $\O$, which sits in $SO(3)$. As we have
seen in Section 3.1, the group $\O=\{[o,o] \mid o\in{\bf O}\}$,
where ${\bf O}={\bf T}\cup\frac{1}{\sqrt2}({\bf 1}+{\bf i}){\bf
T}$ is the octahedral binary subgroup of $U$, of order 48. Since
$\O$ contains $\T$, it follows that $\O$ is generated by the
generators of $\T$ and the extra element $\l[\frac{1}{\sqrt2}({\bf
1}+{\bf i}),\frac{1}{\sqrt2}({\bf 1}+{\bf i})\r]$ with
representing matrix
$\biggl(\begin{smallmatrix}1&0&0&0\\0&1&0&0\\0&0&0&1\\0&0&-1&0
 \end{smallmatrix}\biggl)$. This extra element fixes the
polynomial (\ref{polyhedrap}) if and only if $s=0$. Therefore, the
harmonic cubics stabilized by $\O$ lie on the $\so$-orbit of the
polynomial $rx_1(x_1^2-x_2^2-x_3^2-x_4^2), r>0$  which has full
symmetry $SO(3)$.

We look now for harmonic cubics invariant under the irreducibly
acting octahedral subgroup $\O^+=\{[o,o],o\in{\bf T} \mbox{
  and } [o,-o],o\in \frac{1}{\sqrt2}({\bf 1}+{\bf i}){\bf T}\}$.
The subgroup $\O^+$ contains $\T$ and it is generated by the
generators of $\T$ plus the extra element
$\l[\frac{1}{\sqrt2}({\bf 1}+{\bf i}),-\frac{1}{\sqrt2}({\bf
1}+{\bf i})\r]$. This extra element fixes the harmonic polynomial
(\ref{polyhedrap}) if and only if $r=0$. Therefore, the harmonic
cubics stabilized by $\O^+$ lie on the $\so$-orbit of the
polynomial $sx_2x_3x_4, \ s>0$.

We look now for the harmonic cubics invariant under the reducibly
acting icosahedral subgroup $\I$, which sits in $SO(3)$. As we
have seen in Section 3.1, $\I=\{[l,l]\mid l\in{\bf I}\}$, where
${\bf I}=\cup_{k=0}^{4}(\frac{1}{2\tau}+\frac{\tau}{2}{\bf
i}+\frac{1}{2}{\bf j})^k{\bf T}$ is the binary icosahedral
subgroup of $U$, of order 120 and $\tau=\frac{\sqrt5 +1}{2}$. The
subgroup $\I$ contains $\T$ and it is generated by the generators
of $\T$ plus the extra element
$\l[\frac{1}{2\tau}+\frac{\tau}{2}{\bf i}+\frac{1}{2}{\bf
j},\frac{1}{2\tau}+\frac{\tau}{2}{\bf i}+\frac{1}{2}{\bf j}\r]$.
Straightforward calculations show that this extra element fixes
the harmonic polynomial (\ref{polyhedrap}) if and only if $s=0$.
Therefore, the harmonic cubics stabilized by $\I$ lie on the
$\so$-orbit of the polynomial $rx_1(x_1^2-x_2^2-x_3^2-x_4^2)$
which has full symmetry $SO(3)$.

Finally, we look for the harmonic cubics invariant under the
irreducibly acting icosahedral subgroup $\I^+$. From Section 3.1,
$\I^+=\{[r^+,r] \mid r\in{\bf I}\}$, where $r^+$ is the image of
$r\in{\bf I}$ under the automorphism of the quaternion field that
changes the sign of $\sqrt 5$. This automorphism exchanges $\tau$
for $-\frac{1}{\tau}$. The subgroup $\I^+$ contains $\T$ and it is
generated by the generators of $\T$ plus the extra element
$\l[(\frac{1}{2\tau}+\frac{\tau}{2}{\bf i}+\frac{1}{2}{\bf
j})^+,\frac{1}{2\tau}+\frac{\tau}{2}{\bf i}+\frac{1}{2}{\bf j} \
\r]=\l[-\frac{\tau}{2}-\frac{1}{2\tau}{\bf i}+\frac{1}{2}{\bf
j},\frac{1}{2\tau}+\frac{\tau}{2}{\bf i}+\frac{1}{2}{\bf j}\r]$.
Straightforward calculations show that this extra element fixes
the harmonic polynomial (\ref{polyhedrap}) if and only if
$s=2\sqrt 5 r$. Therefore, the harmonic cubics stabilized by
$\I^+$ lie on the $\so$-orbit of the polynomial
$r[x_1(x_1^2-x_2^2-x_3^2-x_4^2)+2\sqrt5 x_2x_3x_4]$.

To conclude, one can easily compute that the identity component of
the stabilizer of the polynomial (\ref{polyhedrap}) is always
discrete, except in the case $s=0$.  $\Box$

\vskip.4cm

We now consider those special Lagrangian submanifolds $L\subset
\C^4$ whose cubic fundamental form has a polyhedral symmetry at
each point.
\begin{theorem}
Suppose that $L\subset \C^4$ is a connected special Lagrangian
$4$-fold with the property that its fundamental cubic at each
point has a tetrahedral symmetry $\T$. Then, up to congruence, $L$
is the Harvey-Lawson example $L\subset \C^4$ defined in standard
coordinates by the equations
\begin{align*}
L: \quad &|z_0|=|z_1|=|z_2|=|z_3| \\
&\mbox{Re}(z_0z_1z_2z_3)=\sqrt[5]{2}
\end{align*}
\end{theorem}

\pf Let $L$ be a special Lagrangian 4-fold that satisfies the
hypotheses of the theorem and let $C$ be its fundamental cubic.
Proposition \ref{polyhedral} implies that there exist functions
$r,s\colon L\to \R_+$ with $s\not= \sqrt 5 r$ and a $\T$-subbundle
$F\subset P_L$ over $L$ for which the following identity holds:
$$ C=3r\om_1(\om_1^2-\om_2^2-\om_3^2-\om_4^2-)+6s\om_2\om_3\om_4 $$

Since $F$ is an $\T$-bundle, the following relations hold:
$\a_{21}=\a_{31}=\a_{41}=\a_{32}=\a_{42}=\a_{43}=0 \mod
\{\om_1,\om_2,\om_3,\om_4\}$. The usual differential analysis
yields the following structure equations on $F$:
\begin{align}
d\om_1&=0\notag\\
d\om_2&=\sqrt{s^2-r^2}\om_1\wedge\om_2 \notag\\
d\om_3&=\sqrt{s^2-r^2}\om_1\wedge\om_3\label{eq41}\\
d\om_4&=\sqrt{s^2-r^2}\om_1\wedge\om_4,\notag\\
dr&=-5r\sqrt{s^2-r^2}\om_1,\notag\\
ds&=-s\sqrt{s^2-r^2}\om_1 \notag
\end{align}

From the last two equations in (\ref{eq41}), it follows that
$r=cs^5$, $c>0$ constant. We can suppose that $c=1$ since the
equations are invariant under scaling. Moreover, $s\in[-1,1]$.

The above structure equations imply that $\om_1=0$ defines an
integrable 3-plane field which we denote by
 $\Gamma_2$ and that $\om_2=\om_3=\om_4=0$ defines an integrable 1-plane field denoted by $\Gamma_1$. Since $d\om_1=0$,
 it follows that $\om_1=dx_1$ on the leaves of the
 foliation $\Gamma_1$. The structure equations (\ref{eq41}) also imply $d(s\om_2)=d(s\om_3)=d(s\om_4)=0$ and therefore
  there exist functions $x_2,x_3,x_4$ on $L$ such that $\om_2=\frac{dx_2}{s}, \ \om_3=\frac{dx_3}{s}$ and
  $\om_4=\frac{dx_4}{s}$. The metric $g=\frac{dx_2^2+dx_3^2+dx_4^2}{s^2}$ is well defined on the leaves
  of the $\Gamma_2$ foliation.

Equations $de_i = e_j\al_{ji}-Je_j\b_{ji}$ and $d(Je_i)=
e_j\b_{ji}+Je_j\al_{ji}$ yield that, as matrices
$$ d(e_1 \quad Je_1)=(e_1 \quad Je_1)\begin{pmatrix}0&3s^5\om_1\\-3s^5\om_1&0\end{pmatrix} \mod\{\om_2,\om_3,\om_4\}$$
Therefore, the leaves of the $\Gamma_1$ foliation are plane curves
with curvature $k=3s^5$, lying in the complex line $(e_1, Je_1)$.
These curves are congruent since $ds$ is a multiple of $\om_1$.

Now consider the $\Gamma_2$ foliation, defined by the equation
$\om_1=0$. Since $s$ is
constant on its leaves, the connection matrix $A=\begin{pmatrix} \a_{ij}&\b_{ij}\\
-\b_{ij}& \a_{ij}\end{pmatrix}$ satisfies $A\wedge A=dA=0$.
Therefore $A$ takes values in a 3-dimensional abelian subalgebra
${\mathfrak g}\subset {\mathfrak su}(4)$. The maximal torus of
$SU(4)$ is conjugate to the subgroup $$T^3=\l\{\mbox
{diag}(e^{i\theta_0},e^{i\theta_1},e^{i\theta_2},e^{i\theta_3})\mid
\sum_{k=0}^{3} \theta_i=0 \mod 2\pi\r\}$$ and the maximal torus
acts on $L$ by rotating around a plane curve $C$. Therefore, the
solution is invariant under the torus action and the only special
Lagrangian 4-folds with this property are described explicitly by
Harvey and Lawson in their paper \cite{hl}. If $(z_0,z_1,z_2,z_3)$
are coordinates on $\C^4$, then the special Lagrangian 4-folds in
$\C^4$ invariant under $T^3$ look like:
\begin{align}
|z_0|^2-|z_1|^2=c_1, \quad &|z_0|^2-|z_2|^2=c_2, \quad |z_0|^2-|z_3|^2=c_3,\label{eq42} \\
&\mbox{Re}(z_0z_1z_2z_3)=a \notag
\end{align}
for some real constants $a,c_1,c_2,c_3$. It is easy to see that
the solution of the structure equations (\ref{eq41}) is symmetric
in $(z_1,z_2,z_3)$, therefore $c_1=c_2=c_3=c$.

Reparametrizing the solution using polar coordinates $z_k=r_k
e^{i\theta_k}$, $k=0...3$, (\ref{eq42}) becomes
\begin{align}
&r_0^2-r_1^2=r_0^2-r_2^2=r_0^2-r_3^2=c \label{eq43} \\
&\theta_0=\arccos\frac{a}{r_0r_1r_2r_3}-\theta_1-\theta_2-\theta_3
\notag
\end{align}

We will find out for what constants $a$ and $c$, the special
Lagrangian 4-fold defined by (\ref{eq43}) is a solution of  the
structure equations. As we have seen, the solution is a special
Lagrangian 4-fold which is foliated by congruent curves of
curvature $3s^5$ and 3-manifolds which are the $T^3$-orbit of the
points on the leaves of the first foliation.

Since $z(r,\theta_1,\theta_2,\theta_3)=(\sqrt{c+r^2},
\arccos\frac{a}{r^3\sqrt{c+r^2}}-\theta_1-\theta_2-\theta_3,r,\theta_1,r,\theta_2,r,\theta_3)$,
the tangent plane to a $T^3$-orbits is spanned by the vectors
$v_1=z_{\theta_1}=-\frac{\partial}{\partial
\theta_0}+\frac{\partial}{\partial \theta_1}, \
v_2=z_{\theta_2}=-\frac{\partial}{\partial
\theta_0}+\frac{\partial}{\partial \theta_2}$ and $
v_3=z_{\theta_3}=-\frac{\partial}{\partial
\theta_0}+\frac{\partial}{\partial \theta_3}$.

We look now for another vector $v_0$ in the tangent space of $L$
such that $\{ v_0,v_1,v_2,v_3\}$ are a basis of this tangent
space. Since this tangent space is special Lagrangian, the
symplectic form $\om$ and the imaginary part of the holomorphic
volume form $\Omega$ should vanish on it. Also, $v_0$ should be
orthogonal to $v_i, \ i=1...3$, so $g(v_0,v_i)=0$ for $i=1...3$.

Let us write
\begin{equation*}
v_0=\sum_{i=0}^{3}\mu_i \frac{\partial}{\partial \theta_i}+\nu_i
\frac{\partial}{\partial r_i}
\end{equation*}
The symplectic form in polar coordinates is
$\om=\sum_{i=0}^{3}r_idr_i\wedge d\theta_i$ and the condition
$\om(v_0,v_i)=0$ for $i=1,2,3$ implies that
$\nu_i=\frac{r_0\nu_0}{r_i}=\frac{\nu}{r_i}$, where
$\nu=r_0\nu_0$. The metric on $\C^4$ is
$g=\sum_{i=0}^{3}(dr_i)^2+r_i^2(d\theta_i)^2$ and the condition
that $g(v_0,v_i)=0$ for $i=1...3$ implies the relations
$\mu_i=\frac{\mu_0r_0^2}{r_i^2}=\frac{\mu}{r_i^2}$, for $i=1,2,3$,
where $\mu=\mu_0r_0^2$. Finally, straightforward calculations
yield that the condition $\Im \Omega(v_1,v_2,v_3,v_0)=0$ implies
the relation
$\frac{\nu}{\mu}=\tan(\theta_0+\theta_1+\theta_2+\theta_3)$.Therefore,
\begin{equation*}
v_0=\sum_{i=0}^{3}\frac{1}{r_i^2}\frac{\partial}{\partial
\theta_i}+\frac{\tan(\theta_0+\theta_1+\theta_2+\theta_3)}{r_i}\frac{\partial}{\partial
r_i}
\end{equation*}

Next, we find an integral curve of the vector field $v_0$, that
lies in:
\begin{align*}
 L:  \quad  &r_1=r_2=r_3=r, \
 r_0=\sqrt{c+r^2}, \ \theta_0=\arccos\frac{a}{r^3\sqrt{c+r^2}}-\theta_1-\theta_2-\theta_3.
\end{align*}
 When $c=0$, an
integral curve is given by
\begin{align*}
C: \quad &r_0=r_1=r_2=r_3=r, \
\theta_0=\theta_1=\theta_2=\theta_3=\theta, \
\cos(4\theta)=\frac{a}{r^4}
\end{align*}
Therefore, the curve $C$ is given by:
$z_1=z_2=z_3=z_4=re^{i\theta}$, $r^4\cos(4\theta)=a$ and it is a
plane curve which lies in the complex line $z_1=z_2=z_3=z_4$. We
have seen that $L$ is foliated by plane curves with curvature
$3s^5$, so we will determine for what value of $a$ the curve $C$
has this curvature. We choose an orthonormal basis in the
$z_1=z_2=z_3=z_4$ plane: ${\bf
e_1}=(\frac{1}{2},0,\frac{1}{2},0,\frac{1}{2},0,\frac{1}{2},0)$
and ${\bf
e_2}=(0,\frac{1}{2},0,\frac{1}{2},0,\frac{1}{2},0,\frac{1}{2})$
and in this basis, the curve $C$ is given by
$\gamma(\theta)=(2r\cos\theta,2r\sin\theta)$, where
$r=\bigl(\frac{a}{\cos(4\theta)}\bigl)^{\frac{1}{4}}$.

 Computing the
curvature of $\gamma$, one gets
$k(\theta)=-\frac{3}{2}a^{-\frac{1}{4}}(\cos(4\theta))^{\frac{5}{4}}$.
But the curves that foliate $L$ are parameterized by arclength
\begin{align}
&\om_1=dt=d\theta|\gamma \ '|=2a^{\frac{1}{4}}(\cos
 \theta)^{-\frac{5}{4}}\label{eq44}
\end{align}
and the curvature in this parameterization is $k=3s^5$. Therefore
\begin{equation}
\label{eq45}
s=(\frac{k}{3})^{\frac{1}{5}}=-\frac{1}{2^{\frac{1}{5}}a^{\frac{1}{20}}}(\cos4\theta)^{\frac{1}{4}}.
\end{equation}
From the structure equations (\ref{eq41}), it follows that
$ds=-s\sqrt{s^2-s^{10}}\om_1$ has to be satisfied. Using equations
(\ref{eq44}) and (\ref{eq45}), we get
\begin{equation}
\label{eq46} ds=\frac{1}{2^{\frac{6}{5}}a^{\frac{3}{10}}}(\cos
4\theta)^{\frac{1}{2}}\sin 4\theta \ \om_1
\end{equation}
and from equation (\ref{eq46}),
\begin{equation}
\label{eq47}
-s\sqrt{s^2-s^{10}}\om_1=\frac{1}{2^{\frac{1}{5}}a^{\frac{1}{20}}}(\cos
4\theta)^{\frac{1}{4}}\l(\frac{1}{2^{\frac{2}{5}}a^{\frac{1}{10}}}(\cos
4\theta)^{\frac{1}{2}}-\frac{1}{4a^{\frac{1}{2}}}(\cos
4\theta)^{\frac{5}{2}}\r)^{\frac{1}{2}}\om_1
\end{equation}
Equating these last two equations, it follows that
$a=\sqrt[5]{2}$. The structure equations are satisfied now and $L$
is a special Lagrangian 4-fold.

To conclude, the special Lagrangian submanifold $L$ that is a
solution to the structure equations (\ref{eq41}) can be described
explicitly as:
\begin{equation*}
L: \quad |z_0|=|z_1|=|z_2|=|z_3|, \quad
\mbox{Re}(z_0z_1z_2z_3)=\sqrt[5]{2}. \quad \Box
\end{equation*}

\begin{theorem}
Suppose that $L\subset \C^4$ is a connected special Lagrangian
$4$-fold with the property that its fundamental cubic at each
point has an octahedral symmetry $\O^+$, where $\O^+$ is the
irreducibly acting octahedral subgroup of $\so$. Then, up to
congruence, $L$ is the Harvey-Lawson cone in $\C^4$ defined in
standard coordinates $(z_0,z_1,z_2,z_3)$ by the equation
\begin{equation}
\label{harla} L: \quad |z_0|=|z_1|=|z_2|=|z_3|,\quad
\mbox{Re}(z_0z_1z_2z_3)=0
\end{equation}
\end{theorem}

\pf Let $L$ be a special Lagrangian 4-fold that satisfies the
hypotheses of the theorem and let $C$ be its fundamental cubic.
Proposition \ref{polyhedral} implies that there exists a function
$s\colon L\to \R_+$ and an $\O^+$-subbundle $F\subset P_L$ over
$L$ for which the following identity holds:
$$ C=6s\om_2\om_3\om_4 $$ and the 1-forms $\om_1,\om_2,\om_3,\om_4$ form a
basis on $F$.

Straightforward calculations show that the structure equations
are:
\begin{align}
%\label{eq54}
 &d\om_1=0, \ d\om_2=s\om_1\wedge\om_2, \
d\om_3=s\om_1\wedge\om_3, \ d\om_4=s\om_1\wedge\om_4, \
ds=-s^2\om_1\label{eq54}
\end{align}
The structure equations imply the equation
$$de_1=s(e_2\om_2+e_3\om_3+e_4\om_4)=s(dx-e_1\om_1),$$ where $x\colon L^+\to \C^4$. From
here and the last equation in (\ref{eq54}), it follows that
$x=\frac{e_1}{s}+x_0$, where $x_0$ is a constant which we can
reduce to 0 by translation. Therefore $x=\frac{e_1}{s}$. On the
leaves of the foliation $\om_2=\om_3=\om_4=0$,  $de_1=0$ and thus
the vector $e_1$ is constant along these leaves. This tells us
that the special Lagrangian 4-fold $L^+$ is a cone on some
3-dimensional manifold $\Sigma\subset S^7$. We have to determine
now for what 3-dimensional manifolds $\Sigma\subset S^7$, the cone
$C(\Sigma)$ is special Lagrangian and satisfies the structure
equations.

In the case $t=-s$ we obtain $x=-\frac{e_1}{s}$ and the solution
is again a cone through the origin, call it $L^-$. We have that
$L=L^+\cup L^-$.

The connection matrix $A=\begin{pmatrix} \a_{ij}&\b_{ij}\\
-\b_{ij}& \a_{ij}\end{pmatrix}$ satisfies $A\wedge A=dA=0$.
Therefore $A$ takes values in an abelian subalgebra ${\mathfrak
g}\subset {\mathfrak su}(4)$. The group $G=\exp{\mathfrak g}$ is a
maximal torus of $SU(4)$ and it is conjugate to the diagonal torus
$T^3=\l\{\mbox
{diag}(e^{i\theta_0},e^{i\theta_1},e^{i\theta_2},e^{i\theta_3}):
\sum_{k=0}^{3} \theta_i=0 \mod 2\pi\r\}$.

 $G$ acts transitively on
the cone, so the cone is homogeneous and we have to determine
which of the orbits on the 7-sphere are special Lagrangian. The
solution is invariant under the torus action and therefore the
links of these special Lagrangian cones are 3-dimensional tori on
$S^7$. They are described explicitly by Harvey and Lawson in their
paper \cite{hl}. It follows that $L$ is given in standard
coordinates $(z_0,z_1,z_2,z_3)$ by (\ref{harla}).

Therefore, the special Lagrangian cone $L$ is a union of two cones
$L^+$ (obtained in the case $t=s$) and $L^-$ (obtained in the case
$t=-s$) with vertices at the origin through the 3-dimensional tori
$T^+$ and $T^-$ on $S^7$ given by
\begin{align*}
&T^+=\l\{\l(\frac{1}{2}e^{i\theta_0},\frac{1}{2}e^{i\theta_1},\frac{1}{2}e^{i\theta_2},\frac{1}{2}e^{i\theta_3}\r):
\theta_0+\theta_1+\theta_2+\theta_3=\frac{\pi}{2}\r\}\\
&T^-=\l\{\l(\frac{1}{2}e^{i\theta_0},\frac{1}{2}e^{i\theta_1},\frac{1}{2}e^{i\theta_2},\frac{1}{2}e^{i\theta_3}\r):
\theta_0+\theta_1+\theta_2+\theta_3=\frac{3\pi}{2}\r\}
\end{align*}

\begin{theorem}
There are no connected special Lagrangian $4$-folds whose
fundamental cubic at each point has an icosahedral symmetry
$\I^+$, where $\I^+$ is the irreducibly acting icosahedral
subgroup of $\so$.
\end{theorem}

\pf Let $L$ be a special Lagrangian 4-fold that satisfies the
hypotheses of the theorem and let $C$ be its fundamental cubic.
Proposition \ref{polyhedral} implies that there exists a function
$r\colon L\to \R_+$ for which the equation
$$ C= 3r[\om_1(\om_1^2-\om_2^2-\om_3^2-\om_4^2)+2\sqrt5 \om_2\om_3\om_4]$$
defines an $\I^+$-subbundle $F\subset P_L$ of the $L$-adapted
coframe bundle $P_L \ra L$. The usual differential analysis on the
subbundle $F$ yields $r=0$, contrary to the hypothesis. $\Box$

\subsubsection{Symmetries of order $6$, $5$ and $4$}

We have seen in Corollary \ref{cor1} that the elements of a
discrete stabilizer of a fundamental cubic of a special Lagrangian
4-fold have order less or equal to 6. From the proof of this
corollary, the general harmonic cubic stabilized by an element of
order 6 is
$$C=\Re(rz_1^3+sz_1z_2^2)=r(x_1^3-3x_1x_2^2)+s[(x_3^2-x_4^2)x_1-2x_2x_3x_4]$$
where we can arrange $r,s$ to be real and non-negative by making
rotations in the $z_1$ and in the $z_2$-lines. Easy computations
show that the full stabilizer of $C$ is the dihedral group on 6
elements ${\bf D_6}$, if $r\not= 0$ and $s\not= 0$. A similar
differential analysis as in the previous cases yields the
following result.
\begin{theorem}
\label{theoremz6}
 There are no nontrivial special Lagrangian submanifolds in $\C^4$
 whose fundamental cubic has a discrete stabilizer which contains at
 least an element of order $6$.
\end{theorem}

Next, the general harmonic cubic stabilized by an element of order
5 is $$C=\Re(rz_1z_2^2+sz_1^2{\bar z_2})$$ where we can arrange
$r,s\geq 0$. Same analysis gives:
\begin{theorem}
 There are no nontrivial special Lagrangian submanifolds in $\C^4$
 whose fundamental cubic has a discrete stabilizer which contains at
 least an element of order $5$.
\end{theorem}
{\it Remark:} From this theorem, the result in Theorem 3.3.6
follows immediately, since the irreducibly acting icosahedral
subgroup of $\so$ has elements of order 5.

Next, the general harmonic cubic stabilized by an element of order
4 is
$$C=\Re(r{\bar z_1}z_2^2+sz_1z_2^2)$$ where we can arrange again
$r,s$ to be real and non-negative. The stabilizer of $C$ is a
continuous subgroup if $r=0$ or $s=0$, the irreducibly acting
octahedral subgroup $\O^+$ if $r=s$ and the dihedral group ${\bf
D_4}$ in the rest of the cases, since the element of order 2 that
flips the signs of $\{x_2,x_3\}$ belongs to the stabilizer.

We obtain:
\begin{theorem}
\label{theoremz4} There is no nontrivial special Lagrangian
$4$-fold in $\C^4$ whose fundamental cubic has a ${\bf
D_4}$-symmetry at each point.
\end{theorem}

For the details of the calculations in the above results see
\cite{thesis}.

\subsubsection{Discrete symmetry at least $\Z_3$}

Now we consider those special Lagrangian 4-folds $L\subset \C^4$
whose fundamental cubic has at least a $\Z_3$-symmetry at each
point. We saw in the proof of Corollary \ref{cor1} that there are
two inequivalent orbits that stabilize an element of order 3. We
start with: \vskip.2cm

{\bf Case 1.} $(r,s)=({2 \over 3},0):$ The general harmonic cubic
fixed by the element ${\bf g}=\Bigl(
\begin{smallmatrix}
  e^{\frac{4\pi i}{3}}&0\\ 0& I_2 \end{smallmatrix}\Bigl)$ in the maximal
torus is:
$$C=\Re(rz_1^3+tz_2^3+sz_2(|z_2|^2-2|z_1|^2)),$$
where $r,t,s\in\C$. By rotations in the $z_1$-line and $z_2$-line,
we can arrange that $r,s$ be real and non-negative. By writing
$t=u+iv, \ u,v\in\R$, the cubic C becomes:
$$C=r(x_1^3-3x_1x_2^2)+sx_3(x_3^3+x_4^2-2x_1^2-2x_2^2)+u(x_3^3-3x_3x_4^2)+v(x_4^3-3x_3^2x_4) \ (*)$$
where $r,u,v,s\in\R$ and $r,s\geq 0$.

The next lemma tells us what the full stabilizer of $C$ is.

\begin{lemma}
The full stabilizer of the harmonic cubic polynomial $(*)$ is:
\vskip.2cm

$1)\ \mbox{a continuous subgroup of $\so$, if $r=0$ or
$s=u=v=0$;}$\vskip.2cm

$2)\ \mbox{the dihedral subgroup ${\bf D_3}$ generated by the
order $3$ element ${\bf g}$ and the
 order $2$  element}$\vskip.1cm
 $\ \mbox { that flips the signs of $\{x_2,x_4\}$, if
$r\not =0, v=0$;}$\vskip.2cm

 $3)\ \mbox{the dihedral subgroup
${\bf D_3}$ generated by the order $3$ element ${\bf g}$  and the
order $2$ element}$\vskip.1cm $ \ \mbox { that flips the signs of
$\{x_2,x_4\}$, if $r\not =0, v=0, u=3s$;}$\vskip.2cm

 $4)\ \mbox{the order $18$ normal subgroup of ${\bf D}_3\times {\bf
D}_3$, if $u=v=0$ and $r,s\not =0$;}$\vskip.2cm

$ 5). \ \mbox{the cyclic subgroup ${\bf \Z_3}$ generated by the
order $3$ element ${\bf g}$ if none of the above }$\vskip.2cm \
$\mbox{relations among the parameters $r,s,u,v$ hold.}$
\end{lemma}

\pf We denoted by $G$ be the stabilizer of the polynomial $C$,
where $r,s\geq 0$. A simple computation shows that $G$ is a
continuous subgroup if and only if $r=0$ or $u=v=s=0$. Therefore,
if $r\not= 0$ and  $s^2+u^2+v^2\not = 0$, the stabilizer $G$ is
discrete.

When $s=0$, we can make a rotation in the $(x_3,x_4)$-plane and
suppose also that $v=0$. In this case, the stabilizer of $C$ is
$G$, the order 18 normal subgroup of ${\bf D}_3\ti {\bf D}_3$
described as follows: let the first ${\bf D}_3$ be denoted by
${{\bf D}_3}^+$ and suppose it is generated by the rotation $a_1$
and the reflection $b_1$, where $a_1^3=1,b_1^2=1,\ a_1b_1a_1=b_1.$
Denote the second ${\bf D}_3$ by ${{\bf D}_3}^-$ and suppose it is
generated by the rotation $a_2$ and the reflection $b_2$, where
$a_2^3=1,\ b_2^2=1,\ a_2b_2a_2=b_2.$ Then ${{\bf D}_3}^+$ consists
of the elements
$$\{\theta_1^+=1,\ \theta_2^+=a_1,\ \theta_3^+=a_1^2,\ r_1^+=b_1,\ r_2^+=a_1b_1,\ r_3^+=a_1^2b_1\}$$
and ${{\bf D}_3}^-$ consists of the elements
$$\{\theta_1^-=1,\ \theta_2^-=a_2,\ \theta_3^-=a_2^2,\ r_1^-=b_2,\ r_2^-=a_2b_2,\ r_3^-=a_2^2b_2\}.$$
The $\so$-stabilizer of the cubic $C$ is formed by the 18 pair
elements:$$\{(\theta_i^+,\theta_j^-),(r_i^+,r_j^-),\ i,j=1...3\}$$

Next, if $r\not=0$ and $s\not =0$, the differential analysis
yields the following cases:

i) If $v=0$, the stabilizer $G$ of the cubic $C$ is the dihedral
subgroup ${\bf D_3}$ generated by the order 3 element ${\bf g}$
and the order 2 element that flips the signs of $x_2$ and $x_4$.

ii) If $v=0, u=3s$, the stabilizer $G$ of $C$ is also the above
dihedral subgroup ${\bf D_3}$.

iii) In the general case, when none of the above relations among
the parameters $r,s,u,v$ hold, the stabilizer of $C$ is $\Z_3$.
$\Box$ \smallskip

In the case of ${\bf D_3}$-symmetry we obtain the following
partial result:

\begin{proposition}
 There is an infinite parameter family of connected special
Lagrangian submanifolds in $\C^4$ such that the
 fundamental cubic at each point has a ${\bf D}_3$-symmetry and is of the form
 $(*)$,
 where $v=0$ and $r,s\not =0$. This family depends on $2$ functions of one variable.
\end{proposition}

\pf Let $L$ be a special Lagrangian 4-fold that satisfies the
hypotheses of the theorem and let $C$ be its fundamental cubic. It
follows that
$$C=r(\om_1^3-3\om_1\om_2^2)+u(\om_3^3-3\om_3\om_4^2)+3s\om_3(\om_3^2+\om_4^2-2\om_1^2-2\om_2^2),$$
with $r>0,s>0$ defines a ${\bf {\bf D}_3}$-subbundle $F\subset
P_L$ of the adapted coframe bundle $P_L\to L$. In this case, were
able to write down the structure equations that hold on the bundle
$F$, but were unable to describe completely the family of special
Lagrangian submanifolds in this case. Cartan-K{\"a}hler theorem
tells us that the family should depend on 2 functions of one
variables. For more details see \cite{thesis}. $\Box$

\begin{theorem}
Let $L$ be a connected special Lagrangian submanifolds in $\C^4$
such that its
 fundamental cubic at each point has a ${\bf D}_3$-symmetry and it is of the form $(*)$,
 where $v=0,u=3s$ and $r,s\not =0$. Then $L$ is, up to rigid motion, an open subset
 of the asymptotically conical special Lagrangian $4$-fold given by:
\begin{equation}
\label{assyc}
 L_\Sigma=\{(a+ib){\bf u}|\ {\bf u} \in \Sigma, \ \Re(a+ib)^4=c\},
\end{equation}
  where $c$ is a real
 constant and $\Sigma\subset S^7$ is a
 $3$-manifold with the property that the cone on it is special
 Lagrangian, with phase $i$.
\end{theorem}

\pf Let $L$ be a special Lagrangian 4-fold that satisfies the
hypotheses of the theorem and let $C$ be its fundamental cubic. It
follows that
$$C=r(\om_1^3-3\om_1\om_2^2)+3s\om_3(\om_3^2-\om_1^2-\om_2^2-\om_4^2),$$
with $r>0,s\not= 0$ defines a ${\bf {\bf D}_3}$-subbundle
$F\subset P_L$ of the adapted coframe bundle $P_L\to L$.

The structure equations on this bundle are computed to be:
\begin{align}
d\om_1&= t_1\om_3\wedge\om_1-t_2\om_4\wedge\om_1-2t_5\om_4\wedge\om_2+t_3\om_1\wedge\om_2\notag\\
d\om_2&=t_4\om_1\wedge\om_2-t_1\om_2\wedge\om_3+t_2\om_2\wedge\om_4+2t_5\om_4\wedge\om_1 \notag\\
d\om_3&=0 \notag\\
d\om_4&=6t_5\om_1\wedge\om_2+t_1\om_3\wedge\om_4 \notag\\
dr&=-3rt_4\om_1+3rt_3\om_2-rt_1\om_3+rt_2\om_4 \notag\\
ds&=-5st_1\om_3 \notag\label{eq90}\\
dt_1&=(4s^2-t_1^2)\om_3 \notag \\
dt_2&=m_1\om_1+m_2\om_2-t_1t_2\om_3+(t_1^2+t_2^2+s^2-9t_5^2)\om_4\notag\\
dt_3&=m_3\om_1+(m_4-2r^2+t_1^2+t_2^2+t_3^2+t_4^2+15t_5^2+s^2)\om_2-t_1t_3\om_3+\notag\\
+&\l(t_2t_3-2t_4t_5+\frac{1}{3}m_2\r)\om_4\notag\\
dt_4&=m_4\om_1-(m_3+2t_2t_5)\om_2-t_1t_4\om_3+\l(2t_3t_5+t_2t_4-\frac{1}{3}m_1\r)\om_4\\
dt_5&=\frac{1}{3}m_2\om_1-\frac{1}{3}m_1\om_2-t_1t_5\om_3+2t_2t_5\om_4\notag
\end{align}
for some functions $m_1,m_2,m_3,m_4$. Differentiation of these
equations does not lead to new relations among the quantities. The
differential ideal on the manifold $M=P_L\ti \R^3$ is involutive,
since the Cartan characters can be computed as $s_1=4,
s_2=s_3=s_4=0$ and the space of integral elements at each point is
parameterized by the 4 parameters $m_1,m_2,m_3,m_4$.

The structure equations imply
$d(s^{\frac{8}{5}}+t_1^2s^{-\frac{2}{5}})=0$. Since $F$ and $L$
are connected, it follows that there exists a constant $c>0$ so
that $s^{\frac{8}{5}}+t_1^2s^{-\frac{2}{5}}=c^{\frac{8}{5}}$.
Therefore, there is a function $\theta$, well defined on $L$, that
satisfies:
$$s^{\frac{4}{5}}=c^{\frac{4}{5}}\cos{4\theta}, \quad
s^{-\frac{1}{5}}t_1=c^{\frac{4}{5}}\sin{4\theta}, \quad
|\theta|<\frac{\pi}{8}.$$

From the sixth equation of (\ref{eq90}), it follows that
$$\om_3=\frac{d\theta}{c(\cos{4\theta})^{\frac{5}{4}}}$$
The structure equations imply that $\om_1=\om_2=\om_4=0$ is
integrable and also that $\om_3=0$ defines an integrable 3-plane
field on $L$. The 1-dimensional leaves of the field $\Gamma_1$
defined by $\om_1=\om_2=\om_4=0$ are congruent along $\Gamma_2$,
the codimension 1 foliation defined by $\om_3=0$. This is clear
since:
$$de_3=-3s\om_3Je_3, \quad d(Je_3)=3s\om_3e_3 \ \mod
\{\om_1,\om_2,\om_4\}$$ and $ds=0 \mod \om_3$, meaning $s$ is
constant along each leaf of $\Gamma_2$. The above equations imply
that the leaves of the $\Gamma_1$ foliation are congruent plane
curves of curvature $-3s$, lying in the complex line $(e_3,Je_3)$.

The form of the structure equations tells us that these examples
must be related to the asymptotically conical special Lagrangian
submanifolds, as seen in \cite{br1} for the $\Z_3$-symmetry case
of the special Lagrangian 3-folds.

Suppose that the plane curves which are the leaves of the
$\Gamma_1$ foliation are of the form $\Re
z^{\frac{1}{p}}=c^{\frac{1}{p}}$, where $c$ is a constant and
$p\in\R$ is to be determined. By dilation, we can take $c=1$ and
consider the curve given by $z(t)=(1+it)^p$, in the
$(e_3,Je_3)$-plane. To compute the curvature of this curve, we use
the formula for the curvature in any parametrization and we get:
\begin{equation}
 \label{curv}
 k(t)=\frac{z'\wedge z''}{({\frac{d{\tilde
s}}{dt}})^3e_3\wedge Je_3}=\frac{p-1}{p}(1+t^2)^{-\frac{p+1}{2}}
\end{equation}
Since $k=-3s$, and also using the sixth and the seventh structure
equations in (\ref{eq90}), we compute that $p=\frac{1}{4}$.

Therefore, the leaves of the $\Gamma_1$-foliation are curves given
by the equation $\Re z^4=c$, where $c$ is a constant. From the
equation of the curvature (\ref{curv}), it follows that as $t\ra
\infty$, $k\ra 0$, so these curves flatten out, telling us that
they have an asymptote.

Now we study the $\Gamma_2$-foliation, whose leaves are
3-manifolds. If we set $\theta=0$, i.e. $t_1=0$ and $s=c$, we
obtain a 3-manifold $\Sigma$, immersed in the 7-sphere $S^7$. This
is clear since:
$$d(Je_3)=-s\om_1e_1-s\om_2e_2-s\om_4e_4=-sdx,$$
where $x\colon \Sigma \to \C^4$ is the position vector. Since $s$
is constant on $\Sigma$, it implies that
$$Je_3=-sx+\mbox{ constant},$$ where we can suppose, by translation, that the constant is 0.
Therefore, $x=-\frac{Je_3}{s}$ and $\Sigma$ is immersed in the
7-sphere of radius $\frac{1}{s}$, in the direction $Je_3$.

The structure equations of the leaves of the $\om_3=0$ foliation
are:
\begin{align*}
d\om_1&=-t_2\om_4\wedge\om_1-2t_5\om_4\wedge\om_2+t_3\om_1\wedge\om_2\\
d\om_2&=t_4\om_1\wedge\om_2+t_2\om_2\wedge\om_4+2t_5\om_4\wedge\om_1 \quad \mod {\om_3}\\
d\om_4&=6t_5\om_1\wedge\om_2
\end{align*}
Consider now the following expressions:
\begin{align*}
\eta_i=&s^{\frac{1}{5}}\om_1, \ i=1,2,4,\\
q_i=&s^{-\frac{1}{5}}t_i,\ i=2...5,\\
p=&s^{-\frac{1}{5}}r,\\
v_i=&s^{-\frac{2}{5}}m_i, \ i=1...4.
\end{align*}
The structure equations derived earlier show that
\begin{align}
d\eta_1=&q_2\eta_1\wedge\eta_4+q_3\eta_1\wedge\eta_2+2q_5\eta_2\wedge\eta_4\notag\\
d\eta_2=&q_2\eta_2\wedge\eta_4+q_4\eta_1\wedge\eta_2+2q_5\eta_4\wedge\eta_1\notag\\
d\eta_4=&6q_5\eta_1\wedge\eta_2\notag\\
dq_2=&v_1\eta_1+v_2\eta_2+(\frac{4}{9}+q_2^2-9q_5^2)\eta_4 \label{eq93}\\
dq_3=&v_3\eta_1+(\frac{4}{9}+q_2^2+q_3^2+q_4^2+15q_5^2-2p^2+v_4)\eta_2+\frac{1}{3}(v_2-6q_4q_5+3q_2q_3)\eta_4\notag\\
dq_4=&v_4\eta_1-(v_3+2q_2q_5)\eta_2-\frac{1}{3}(v_1-6q_3q_5-3q_2q_4)\eta_4\notag\\
dq_5=&\frac{1}{3}v_2\eta_1-\frac{1}{3}v_1\eta_2+2q_2q_5\eta_4\notag\\
dp=&-p(3q_4\eta_1-3q_3\eta_2-q_2\eta_4)\notag
\end{align}
Therefore, the metric $g=\eta_1^2+\eta_2^2+\eta_4^2$ is well
defined on each leaf of the $\Gamma_2$-foliation. The
$\theta$-curves meet the 3-manifold $\Sigma$ orthogonally, so it
is easy to see that the image of
$(-\frac{\pi}{8},\frac{\pi}{8})\times \Sigma$ is of the form :
\begin{equation}
\label{as}
 L_{\Sigma}=\{z{\bf u}| \ {\bf u}\in \Sigma, \ z\in\C, \
\Re z^4=c\},
\end{equation}
 where $c$ is a real constant. In order for this to be a
special Lagrangian 4-fold, the cone on the image of $\Sigma$
should be a special Lagrangian 4-fold.

We shall show now that, indeed, the cone on $\Sigma$ is special
Lagrangian with phase $i$. The cone on $\Sigma$ is parameterized
by:
$$(r,z)\ra rz, \ r\in\R^+, \ z\in\Sigma^3.$$

The tangent space to $C(\Sigma)$ has a basis formed by the
vectors:
$$\l(e_1=\frac{\partial}{\partial x_1}, \ e_2=\frac{\partial}{\partial
x_2},\ e_4=\frac{\partial}{\partial x_4},\
Je_3=\frac{\partial}{\partial y_3}\r)$$ Since
$$\om=dx_1\wedge dy_1+dx_2\wedge dy_2+dx_3\wedge dy_3+dx_4\wedge dy_4,$$ it is clear that
$\om\mid_{C(\Sigma)}=0$, so the cone is Lagrangian. Also,
$\Omega=dz_1\wedge dz_2\wedge dz_3\wedge dz_4$ and we can easily
compute that
$$\mbox {Im} \ \Omega\mid _{C(\Sigma)}=dx_1\wedge dx_2\wedge dy_3 \wedge dx_4,$$ which
represents the volume form on the cone, and ${\mbox Re} \
\Omega\mid _{C(\Sigma)}=0$. Therefore, $C(\Sigma)$ is special
Lagrangian with phase $i$. Then, it is well-known \cite{hl} that
(\ref{as}) is a special Lagrangian 4-fold.

\begin{theorem}
Let $L$ be a connected special Lagrangian submanifolds in $\C^4$
such that its
 fundamental cubic at each point has a $G$-symmetry, where $G$ is the order $18$ normal subgroup of
 ${\bf D}_3\times {\bf D}_3$. Then $L$ is congruent to the product of two holomorphic curves in $\C^2$.
\end{theorem}

\pf Let $L$ be a special Lagrangian 4-fold that satisfies the
hypotheses of the theorem and let $C$ be its fundamental cubic.
Then
$$C=r(\om_1^3-3\om_1\om_2^2)+v(\om_3^3-3\om_3\om_4^2),$$
with $r>0,v>0,r\not =v$ defines a $G$-subbundle $F\subset P_L$ of
the adapted coframe bundle $P_L\to L$, where $G$ is the order 18
normal subgroup of ${\bf D}_3\times {\bf D}_3$.

The structure equations on the subbundle $F$ are computed to be:
\begin{align}
d\om_1&= t_3\om_1\wedge\om_2\notag, \ d\om_2=t_4\om_1\wedge\om_2\notag\\
d\om_3&=t_1\om_3\wedge\om_4, \ d\om_4=t_2\om_3\wedge\om_4 \notag\\
dr&=-3rt_4\om_1+3rt_3\om_2 \notag\\
dv&=-3vt_2\om_3+3vt_1\om_4 \label{eq99}\\
dt_1&=u_1\om_3+(t_1^2+t_2^2-2v^2+u_2)\om_4\notag \\
dt_2&=u_2\om_3-u_1\om_4\notag\\
dt_3&=u_3\om_1+(t_3^2+t_4^2-2r^2+u_4)\om_2\notag\\
dt_4&=u_4\om_1-u_3\om_2\notag
\end{align}
for some functions $u_1,u_2,u_3,u_4$.

From the above structure equations, we can see that
$\om_1=\om_2=0$ and $\om_3=\om_4=0$ define integrable 2-plane
fields on $L$. The structure equations also show that the leaves
of the 2-plane field $\Gamma_1$ defined by $\om_3=\om_4=0$ are
congruent along $\Gamma_2$, the codimension 2 foliation defined by
$\om_1=\om_2=0$. Also, the 2-dimensional leaves of the 2-plane
field $\Gamma_2$ are congruent along $\Gamma_1$.

Since $d(e_1\om_1+e_2\om_2)=0$ and $d(e_3\om_3+e_4\om_4)=0$, it
follows that $e_1\om_1+e_2\om_2=d\pi_1$ and
$e_3\om_3+e_4\om_4=d\pi_2$, where the projections $\pi_1\colon
L\to  \Sigma_1$ and $\pi_2\colon L\to  \Sigma_2$ are well defined.
Therefore, $x=\pi_1+\pi_2+\mbox {const}$ and $L$ is the sum of two
surfaces: $ L=\Sigma_1\times \Sigma_2$.

 Since $d(e_1\wedge e_2\wedge Je_1\wedge Je_2)=0$, it follows that $\Sigma_1$
lies in the complex plane $(e_1,e_2,Je_1,Je_2)$. Also, $\Sigma_2$
lies in the complex plane $(e_3,e_4,Je_3,Je_4)$.

Because $L$ is special Lagrangian, both surfaces $\Sigma_1$ and
$\Sigma_2$ should be special Lagrangian 2-folds in $\C^2$. It is
well-known then that these surfaces should be holomorphic curves
with respect to some complex structure on $\C^2$. More explicitly,
if $\Sigma_1\subset \C^2$, with complex coordinates
$\{z_1=x_1+iy_1,z_2=x_2+iy_2\}$, then $\Sigma_1$ is a holomorphic
curve with respect to the complex coordinates
$\{u_1=x_1-ix_2,v_1=y_1+iy_2\}$. If $\Sigma_2\subset \C^2$, with
standard complex coordinates $\{z_3=x_3+iy_3,z_4=x_4+iy_4\}$, then
$\Sigma_2$ is a holomorphic curve with respect to the complex
coordinates $\{u_2=x_3-ix_4,v_2=y_3+iy_4\}$.    $\Box$
 \smallskip

In the next case of $\Z_3$-symmetry we were unable to describe
completely the SL 4-folds and therefore we have only a partial
result.
\begin{proposition}
There is an infinite parameter family of connected special
Lagrangian submanifolds in
  $\C^4$ such that the
 fundamental cubic at each point has a ${\bf Z}_3$-symmetry and is of the form $(*)$.
  The family depends on $4$ functions of one variable and the elements of this
   family are foliated by non-congruent minimal Legendrian
  surfaces in the direction $\{\om_1,\om_2\}$ and by congruent holomorphic
  curves in the direction $\{\om_3,\om_4\}$.
\end{proposition}

\pf Let $L$ be a special Lagrangian 4-fold that satisfies the
hypotheses of the theorem and let $C$ be its fundamental cubic.
Then
$$C=r(\om_1^3-3\om_1\om_2^2)+3s(\om_3^2+\om_4^2-2\om_1^2-2\om_2^2)\om_3+
u(\om_3^3-3\om_3\om_4^2)+v(\om_4^3-3\om_3^2\om_4)$$ with $r>0$
defines a ${\bf {\bf Z}_3}$-subbundle $F\subset P_L$ of the
adapted coframe bundle $P_L\to L$. The differential analysis shows
that the structure equations on the bundle $F$ are:
\begin{align}
d\om_1&= t_3\om_1\wedge\om_2,\ d\om_2=t_4\om_1\wedge\om_2,\
d\om_3=t_1\om_3\wedge\om_4,\ d\om_4=t_2\om_3\wedge\om_4, \notag\\
dr&=-3rt_4\om_1+3rt_3\om_2+2rst_2\om_3+2rst_1\om_4 \notag\\
ds&=[-svt_1+s(7s+u)t_2]\om_3+s[(3s-u)t_1-vt_2]\om_4\notag\\
du&=[-11svt_1+(3v^2+3u^2+6s^2-su)t_2+t_6]\om_3+\notag \\
[&(11su-3u^2-6s^2-3v^2)t_1-svt_2-t_5]\om_4\notag\\
dv&=t_5\om_3+t_6\om_4 \notag\\
dt_1&=[v(t_1^2+t_2^2+1)-8st_1t_2]\om_3+[(u-s)(t_1^2+t_2^2+1)]\om_4\label{eq100}\\
dt_2&=-[(u+5s)(t_1^2+t_2^2+1)-8st_1^2]\om_3+v(t_1^2+t_2^2+1)\om_4\notag\\
dt_3&=-m_2\om_1+[4s^2(t_1^2+t_2^2+1)+t_3^2+t_4^2-2r^2+m_1)\om_2+2st_2t_3\om_3+2st_1t_3\om_4\notag\\
dt_4&=m_1\om_1+m_2\om_2+2st_2t_4\om_3+2st_1t_4\om_4\notag\\
dt_5&=[-m_4+3t_1^2(30suv-47s^2v-3u^2v-3v^3)+3t_2^2(-7s^2v+10suv-3v^3-3u^2v)+\notag\\
(&25s-7u)t_1t_6+60sv^2t_1t_2+(7u-3s)t_2t_5-7vt_2t_6-7vt_1t_5-18s^2v-6u^2v+\notag\\
2&4suv-6v^3]\om_3+[m_3+(3s-u)t_1t_5-vt_2t_5+(s-u)t_2t_6+vt_1t_6]\om_4\notag\\
dt_6&=m_3\om_3+m_4\om_4\notag
\end{align}
for some functions $m_1,m_2,m_3,m_4$.

The Cartan-K{\"a}hler analysis tells us that the solution should
depend on 4 functions of 1 variable. From the structure equations,
we can see that $\om_3=\om_4=0$ and $\om_1=\om_2=0$ define
integrable 2-plane fields on $L$. Let $\Gamma_1$ be the
$\om_1=\om_2=0$ foliation and $\Gamma_2$ be the $\om_3=\om_4=0$
foliation. The structure equations of the foliation $\Gamma_1$
show that the leaves are congruent and that the metric
$g_1=\om_3^2+\om_4^2$ is well defined on the leaf space of the
$\Gamma_1$ foliation. It is easy to see that the leaves of the
$\Gamma_2$ foliation are non-congruent.

Notice that if we denote $\Delta^2=4s^2(t_1^2+t_2^2+1)$, then we
get that $\frac{d\Delta}{\Delta}=2s(t_2\om_3+t_1\om_4)$. We see
that $\Delta$ is constant on each leaf of the $\Gamma_2$
foliation. We compute that:
\begin{align}
d(\Delta\om_1)&=t_3\om_1\wedge(\Delta\om_2)=(t_3\om_1+t_4\om_2)\wedge(\Delta\om_2)\label{eq108}\\
d(\Delta\om_2)&=-t_4\om_2\wedge(\Delta\om_1)=-(t_3\om_1+t_4\om_2)\wedge(\Delta\om_1)\notag
\end{align}
and the metric $g_2=(\Delta\om_1)^2+(\Delta\om_2)^2$ is well
defined on the leaf space of the $\Gamma_2$ foliation.

Computations also show that
$d(r^{\frac{1}{3}}\Delta^{\frac{2}{3}}\om_1)=0$ and
$d(r^{\frac{1}{3}}\Delta^{\frac{2}{3}}\om_2)=0$. These imply that
there are functions $x_1,x_2$ on $L$ such that
$r^{\frac{1}{3}}\Delta^{\frac{2}{3}}\om_1=dx_1$ and
$r^{\frac{1}{3}}\Delta^{\frac{2}{3}}\om_2=dx_2$. This gives
$x_1,x_2$ up to additive constants and $x_1+ix_2$ is holomorphic
with respect to the complex structure that $\{\om_1,\om_2\}$
define on the leaf space of $\Gamma_2$.

The above imply that the metric
$$g_2=(\frac{\Delta}{r})^{\frac{2}{3}}(dx_1^2+dx_2^2)=F(x_1,x_2)(dx_1^2+dx_2^2),$$ on the
$\Gamma_2$ leaf space has Gauss curvature:
$$k=1-2(\frac{r}{\Delta})^2=1-\frac{2}{F^3(x_1,x_2)}$$
We obtain a differential equation for the function $F(x_1,x_2)$,
given by:
\begin{equation}
\label{eq109}
 \frac{1}{2}\Delta(\ln F)=\frac{2}{F^2}-F
 \end{equation}
We used here the fact that a metric $ds^2=e^{2u}(dx^2+dy^2)$ is
computed to have the Gauss curvature $K=-\Delta ue^{-2u}$. In our
case, take $u=\frac{1}{2}\ln F$ and (\ref{eq109}) follows. The
function $u$ satisfies the differential equation $\Delta
u=e^{2u}-2e^{-4u}$ (Tzitzeica equation), which is completely
integrable by means of inverse scattering method \cite{sha}. This
is the differential equation satisfied by the curvature of the
metric of a minimal Legendrian immersion in $S^5(1)$, invariant
under $S^1$-action, as shown by Mark Haskin in \cite{ha}, p.14.
Sharipov \cite{sha} shows that the minimal immersion satisfying
Tzitzeica equation are minimal tori which are complexly normal in
$S^5$. Therefore, $L$ is foliated by non-congruent minimal
Legendrian surfaces in the direction $\{\om_1,\om_2\}$ and by
congruent holomorphic curves in the direction $\{\om_3,\om_4\}$.
 We do not have a complete description of the family yet.
\smallskip

We move now to analyze the other orbit that stabilizes an element
of order 3.
\\

 {\bf Case 2.} $(r,s)=(\frac{2}{3},\frac{1}{3})$: This case is equivalent to the $(r,s)=(\frac{1}{3},\frac{1}{3})$
case, when the element in the maximal torus that stabilizes $C$ is
$$g=\Bigl(
\begin{smallmatrix}
  e^{\frac{2\pi i}{3}}&0\\ 0& e^{\frac{2\pi i}{3}}
  \end{smallmatrix}\Bigl),$$ of {\bf order 3}. The general harmonic cubic fixed by this element is
$$C=\Re(a_3z_1^3+3a_2z_1^2z_2+3a_1z_1z_2^2+a_0z_2^3), \
a_0,a_1,a_2,a_3\in\C$$
 We notice that the commutator of $g$ is
larger than the maximal torus in this case. The unitary group
$U(2)$ commutes with $g$ and therefore we can use also its action
to get rid of certain parameters, more precisely to make $a_0=0$
and $a_1,a_2  \in \R$. So, the general cubic stabilized by $g$
will look like:
$$C=u(x_1^3-3x_1x_2^2)+v(3x_1^2x_2-x_2^3)+r[(x_1^2-x_2^2)x_3-2x_1x_2x_4])+s[(x_3^2-x_4^2)x_1-2x_2x_3x_4] (**),$$
where $u,v,r,s\in\R$.

\begin{lemma}
\label{lemmaz3}
 The full stabilizer of the polynomial given by $(**)$, where
$r,s,u,v \in\R$ is \vskip.2cm

$1)\ \mbox{a continuous subgroup of $\so$, if $r=s=0$ or $u=v=r=0$
or $u=v=s=0$;}$\vskip.2cm

 $2)\ \mbox{the dihedral
subgroup ${\bf D}_6$ generated by the order $6$ element ${\bf
a}=\l(
\begin{smallmatrix}
  e^{\frac{4\pi i}{3}}&0\\ 0& e^{\frac{\pi i}{3}}
  \end{smallmatrix}\r)$ and the element }$\vskip.1cm
$ \ \mbox { of order $2$ that flips the signs of $\{x_2,x_4\}$, if
$r=v=0$;}$\vskip.2cm

$3)\ \mbox{ the dihedral subgroup ${\bf D_3}$ generated by the
order $3$ element ${\bf g}$  and the order $2$ element}$\vskip.1cm
$ \ \mbox {  that flips the signs of $\{x_2,x_4\}$, if
$s=v=0$;}$\vskip.2cm

 $4) \ \mbox{the cyclic subgroup ${\bf \Z_3}$
generated by the order $3$ element ${\bf g}$ if none of the above
relations}$\vskip.1cm
 $\ \mbox{  among the parameters $r,s,u,v$
hold.}$
\end{lemma}

\pf We denoted by $G$ be the stabilizer of the polynomial $C$. A
simple computation shows that $G$ is a continuous subgroup if and
only if $r=s=0$ or $u=v=r=0$ or $u=v=s=0$.

Doing the differential analysis in the discrete case, we obtain
the following cases where the stabilizer becomes larger than
$\Z_3$:

i)$r=0$. By making a rotation, if necessary, of angle
$\theta=\frac{1}{3}\arctan(-\frac{v}{u})$ in the $(x_1,x_2)$-plane
and of angle $-2\theta$ in the $(x_3,x_4)$-plane, we can suppose
that $v=0$ also.

The stabilizer of $C$ for $r=v=0$ is seen to be the dihedral
subgroup ${\bf D}_6$ generated by the order 6 element ${\bf
a}=\Bigl(
\begin{smallmatrix}
  e^{\frac{4\pi i}{3}}&0\\ 0& e^{\frac{\pi i}{3}}
  \end{smallmatrix}\Bigl)$ and the element of order 2 that flips the signs of
  $\{x_2,x_4\}$. This case of symmetry at least $\Z_6$ was already studied in
Section 3.3.2 and it did
  not yield any families of special Lagrangian 4-folds.

\vskip.15cm

 ii) If $s=0$. In this case, we can arrange that $v=0$ also and
 the stabilizer is computed to be the dihedral group ${\bf D_3}$.

\vskip.15cm

 iii) In the general case, when none of the above
relations among the parameters $r,s,u,v$ hold, the stabilizer of
$C$ is computed to be $\Z_3$ generated by the element $g$. $\Box$

\begin{theorem}
Let $L$ be a connected special Lagrangian submanifold in $\C^4$
such that its
 fundamental cubic at each point has a ${\bf Z}_3$-symmetry and it is of the form $(**)$.
 Then $L$ is an $I$-special Lagrangian $J$-holomorphic surface in $\C^4$, where
 $\{I,J,K\}$ is the hyper-K{\"a}hler stucture on $\C^4$.
 \end{theorem}

\pf Let $L$ be a special Lagrangian 4-fold that satisfies the
hypotheses of the theorem and let $C$ be its fundamental cubic.
From Lemma \ref{lemmaz3}, the equation
$$C=u(\om_1^3-3\om_1\om_2^2)+v(3\om_1^2\om_2-\om_2^3)
+r[(\om_1^2-\om_2^2)\om_3-2\om_1\om_2\om_4])+s[(\om_3^2-\om_4^2)\om_1-2\om_2\om_3\om_4],$$
with $r,s,u,v\in\R$ defines a $\Z_3$-subbundle $F\subset P_L$ of
the adapted coframe bundle $P_L\to L$.

On the subbundle $F$, the following identities hold:
\begin{equation}
\label{eq120}
(\b_{ij})=\begin{pmatrix} u\om_1+v\om_2+r\om_3&v\om_1-u\om_2-r\om_4&r\om_1+s\om_3&-r\om_2-s\om_4\\
v\om_1-u\om_2-r\om_4&-u\om_1-v\om_2-r\om_3&-r\om_2-s\om_4&-r\om_1-s\om_3\\
r\om_1+s\om_3&-r\om_2-s\om_4&s\om_1&-s\om_2\\
-r\om_2-s\om_4&-r\om_1-s\om_3&-s\om_2&-s\om_1 \end{pmatrix}
\end{equation}

The Cartan-K{\"a}hler analysis yields the following relations
between the $\alpha_{ij}$'s:
$$\a_{31}-\a_{42}=0  \quad \mbox{and} \quad \a_{32}+\a_{41}=0$$

We consider the ideal $I_1$, on the coframe bundle, spanned by the
1-forms (\ref{eq120}) and the two 1-forms $\a_{31}-\a_{42}$ and
$\a_{32}+a_{41}$. The independence condition is given by
$\om_1\wedge \om_2\wedge \om_3\wedge \om_4\not =0$ and the tableau
matrix for the structure equations is given by:
$$\begin{pmatrix} \a_1&\a_2&\a_3&\a_4\\\a_2&-\a_1&\a_4&-\a_3\\
\a_3&\a_4&\a_5&\a_6\\
\a_4&-\a_3&\a_6&-\a_5\\\a_5&\a_6&-3s\pi_6&-3s\pi_5\\
\a_6&-\a_5&-3s\pi_5&3s\pi_6 \end{pmatrix}$$ where $$\pi_1=dr,\
\pi_2=ds,\ \pi_3=du,\ \pi_4=dv,\ \pi_5=\a_{41},\ \pi_6=\a_{42},\
\pi_7=\a_{43},\ \pi_8=\a_{21}$$
\begin{align*}
\a_1&=-\pi_3+3r\pi_6+3v\pi_8\\
\a_2&=-\pi_4-3r\pi_5-3u\pi_8\\
\a_3&=-\pi_1+v\pi_5+(2s-u)\pi_6\\
\a_4&=-(2s+u)\pi_5-v\pi_6-r\pi_7-2r\pi_8\\
\a_5&=-2r\pi_2-\pi_6\\
\a_6&=-2r\pi_5-2s\pi_7+s\pi_8\end{align*} From the above tableau,
we compute the reduced Cartan characters as $s_1'=6,s_2'=2,
s_3'=s_4'=0$. The integral elements of the system at each point is
shown to form a space of dimension $10=s_1'+2s_2'+3s_3'+4s_4'$ and
therefore, by Cartan's Test, the system $I_1$ is involutive.

The form of the tableau resembles the tableau for the structure
equations of a complex surface with complex structure given by
$\gamma_1=\om_1+i\om_2$ and $\gamma_2=\om_3+i\om_4$. We will show
that this is actually the case. The characteristic variety of the
ideal is formed by 2 complex lines spanned by
$\{\gamma_1,\gamma_2\}$ and their conjugates.

The first derived system of $I_1$ is generated by the rank 6
Pfaffian system $I_2$ spanned by the six 1-forms:
\begin{align}
\theta_1&=\b_{11}+\b_{22}, \ \theta_2=\b_{33}+\b_{44},\ \theta_3=\b_{41}-\b_{32} \notag\\
\theta_4&=\b_{31}+\b_{42},\
\theta_5=\a_{31}-\a_{42},\theta_6=\a_{32}+\a_{41}\label{derived}
\end{align}
This system is Frobenius and defines a foliation of dimension 10
on the coframe bundle. The integral manifold of our original
system will be a submanifold of the maximal integral manifold of
the derived system $I_2$. Therefore, we will adapt frames and
restrict to the first derived system, looking for integral
manifolds of this system.

We notice that, when restricted to the first derived system, the
connection matrix takes values in the Lie algebra of a
10-dimensional subgroup of $SU(4)$. This subgroup can be shown to
be $Sp(2)$. The system $I_2$ restricts to the $Sp(2)$-coframe
bundle, of dimension 18. The canonical form on this bundle has
components $\xi_i=\om_i+i\eta_i$ and the 1-forms
$$\{\om_i,\eta_i,\b_{11},\b_{33},\b_{21},\b_{31},\b_{41},\b_{43},\a_{21},\a_{31},
\a_{41},\a_{43}\}$$ form a basis for the space of 1-forms on this
coframe bundle $P\cong \C^4\ti Sp(2)$. On the integral manifolds,
$\eta_i=0$ for $i=1...4$.

Now, the symplectic group $Sp(2)$ leaves invariant 3 symplectic
2-forms $\{\zeta_1,\zeta_2,\zeta_3\}$. One of them is the
K{\"a}hler form of the standard complex structure $I$ on $\R^8$.
\begin{align*}
\zeta_1=\frac{i}{2}(\xi_1\wedge {\bar \xi_1}+\xi_2\wedge {\bar
\xi_2}+\xi_3\wedge {\bar \xi_3}+\xi_4\wedge {\bar
\xi_4})=\om_1\wedge \eta_1+\om_2\wedge \eta_2+\omega_3\wedge
\eta_3+\omega_4\wedge \eta_4
\end{align*}
and the other 2-forms are computed to be:
\begin{align*}
\zeta_2&=\om_1\wedge \om_2+\om_3\wedge\om_4-\eta_1\wedge \eta_2-\eta_3\wedge\eta_4\\
\zeta_3&=\om_1\wedge \eta_2-\om_2\wedge \eta_1+\om_3\wedge
\eta_4-\om_4\wedge \eta_3
\end{align*}

Let $I,J,K$ be the complex structures on $\R^8$ corresponding to
left multiplication by the elementary quaternions $i,j$ and $k$.
Then the standard metric $g=\sum_{i=1}^4(\om_i^2+\eta_i^2)$ on
$\R^8$ is K{\"a}hler with respect to each $I,J,K$, with K{\"a}hler
form $\zeta_1,\zeta_2$ and $\zeta_3$, respectively. The forms
$\psi_1=\zeta_2+i\zeta_3,\ \psi_2=\zeta_1+i\zeta_3, \
\psi_3=\zeta_1+i\zeta_2$ are the holomorphic symplectic forms on
$\C^4$, associated to the complex structures $I,J$ and $K$
respectively. The standard complex structure on $\R^8$ is
considered to be $I$, given by the complex 1-forms:
$\om_j+i\eta_j,\ j=1...4$. The 4-forms
$\Omega_i=\frac{1}{2}\psi_i^2$, $i=1...3$ are the holomorphic
volume forms on $\C^4$, associated to the complex structures $I,J$
and $K$, respectively, with $\Omega_1$ being the usual holomorphic
volume form. On the integral manifolds of $I_2$,
$\zeta_1=\zeta_3=0$ and $\zeta_2=\om_1\wedge
\om_2+\om_3\wedge\om_4$ is the K{\"a}hler form for the complex
structure $J$ given by $\gamma_1=\om_1+i\om_2$ and
$\gamma_2=\om_3+i\om_4$.

We will now show that the integral manifold of the ideal generated
by the 2-forms $\zeta_1$ and $\zeta_3$ are complex manifolds with
respect to the complex structure $J$. Let $(z_1,z_2,z_3,z_4)$ be
the complex coordinates on $\R^8$ that are holomorphic for the
complex structure $J$. Then:
$$\zeta_1+i\zeta_3=dz_1\wedge dz_2+dz_3\wedge dz_4$$
Let ${\bf I}$ be the differential ideal generated by the complex
1-form $\psi_2=\zeta_1+i\zeta_3$. We use the Cartan-K{\"a}hler
analysis \cite{br1} to compute the Cartan characters as
$s_1=s_2=2$ and $s_3=s_4=0$. The space of 2-dimensional integral
elements over a point has dimension $6=s_1+2s_2+3s_3+4s_4$ and by
Cartan's Test, the system is involutive. The maximal integral
manifolds of this ideal are given by 2 complex linear equations,
i.e. they are $J$-holomorphic surfaces in $\C^4$.

An integral manifold of the derived system $I_2$ is an integral
manifold of the system $\zeta_1=\zeta_3=0$ and an integral
manifold of the system $\zeta_1=\zeta_3=0$ is an integral manifold
of the derived system. To summarize, the integral manifolds
$\Sigma$ of our original system are $J$-holomorphic surfaces in
$\C^4$. They are $I$-special Lagrangian 4-folds, because
$$\zeta_1\mid_\Sigma=0 \ \mbox{and} \
\frac{1}{2}\Im(\Omega_1^2)\mid_\Sigma=\zeta_2\wedge \zeta_3\mid
_\Sigma=0. \qquad \Box$$

\begin{theorem}
 Let $L$ be a connected special Lagrangian submanifold in $\C^4$ such that its
 fundamental cubic at each point has a ${\bf D}_3$-symmetry and it is of the form $(*)$  with $s=v=0$.
 Then $L$ is a ruled $I$-special Lagrangian $J$-holomorphic surface in $\C^4$.
\end{theorem}

\pf The analysis here is similar to the one in the previous
result. It can be shown that the solutions
 are again $I$-special Lagrangian J-holomorphic
surfaces. Moreover, the structure equations show that the
holomorphic surfaces are foliated by planes in the
$\{e_3,e_4\}$-direction. The conclusion is that the solutions are
ruled $I$-special Lagrangian $J$-holomorphic surfaces. $\Box$
\\

We conclude this paper with the following:
\smallskip

{\bf Open Problem:}

It remains to study the general case when the symmetry of the
fundamental cubic is at least a $\Z_2$. This is the most
complicated case since the space of fixed harmonic cubics involves
a large number of parameters.

%%%%%%%%%%%%%%%%%%%%%%%%%%%%%%%%%%%%%
%%%%%%%%%%%%%%%%%%%%%%% References
%%%%%%%%%%%%%%%%%%%%%%%%%%%%%%%%%%%%%%%%%%

{\it e-mail:} ionelm@@math.mcmaster.ca

\end{document}